\documentclass{article}
\usepackage{graphicx}
\usepackage{amsmath}
\usepackage{amsfonts}
\usepackage{amssymb}
\newtheorem{theorem}{Theorem}[subsection]

\newtheorem{corollary}[theorem]{Corollary}

\newtheorem{definition}[theorem]{Definition}
\newtheorem{example}[theorem]{Example}

\newtheorem{lemma}[theorem]{Lemma}

\newtheorem{proposition}[theorem]{Proposition}
\newtheorem{remark}[theorem]{Remark}

\newenvironment{proof}[1][Proof]{\textbf{#1.} }{\ \rule{0.5em}{0.5em}}

\begin{document}

\title{Schensted type correspondence for type $G_{2}$ and computation of the
canonical basis of a finite dimensional $U_{q}(G_{2})$-module}
\author{C\'{e}dric Lecouvey\\lecouvey@math.unicaen.fr}
\date{}
\maketitle
\begin{abstract}
We use Kang-Misra's combinatorial description of the crystal graphs for
$U_{q}(G_{2})$ to introduce the plactic monoid for type $G_{2}$. Then we
describe the corresponding insertion algorithm which yields a Schensted type
correspondence. Next we give a simple algorithm for computing the canonical
basis of any finite dimensional $U_{q}(G_{2})$-module.
\end{abstract}

\section{Introduction}

The quantum algebra $U_{q}(\frak{g})$ associated to a semisimple Lie algebra
$\frak{g}$ is the $q$-analogue introduced by Drinfeld and Jimbo of its
universal enveloping algebra $U(\frak{g})$. According to Kashiwara \cite{Ka}
each finite dimensional irreducible $U_{q}(\frak{g})$-module $V$ has a unique
crystal basis which can be regarded as a basis at $q=0$.\ This crystal basis
can be extended to obtain a true basis of $V$ called the global crystal basis
of $V.\;$The global crystal basis coincides with the canonical basis of $V$
introduced by Lusztig \cite{Lut}. Moreover the module structure on $V$ induces
a combinatorial structure on its crystal basis $V$ called crystal graph.
Crystal graphs permit to reduce many problems in the representation theory of
$U_{q}(\frak{g})$ to combinatorics.

In this article we restrict ourselves to $U_{q}(G_{2}).\;$Write $\Lambda_{1}$
and $\Lambda_{2}$ for the two fundamentals weights of $U_{q}(G_{2})$ and let
$P_{+}$ be the set of its dominants weights.\ For $\lambda\in P_{+}$ denote by
$V(\lambda)$ the $U_{q}(G_{2})$-module of weight $\lambda$ and by $B(\lambda)$
the crystal graph of $V(\lambda).$ The purpose of this article is two-fold.

In a first part we use the explicit description of the crystal graphs for type
$G_{2}$ \cite{KM} (based on a notion of tableaux for type $G_{2}$) to
introduce a monoid structure $\equiv$ on the vertices of $\Gamma
=\underset{l\geq0}{%
%TCIMACRO{\tbigoplus }%
%BeginExpansion
{\textstyle\bigoplus}
%EndExpansion
}B(\Lambda_{1})^{\otimes l}$ analogous to the plactic monoid of Lascoux and
Sch\"{u}tzenberger.\ Given $w_{1}$ and $w_{2}$ two vertices of $\Gamma,$ this
monoid structure is such that $w_{1}\equiv w_{2}$ if and only if the vertices
$w_{1}$ and $w_{2}$ occur at the same place in two isomorphic connected
components of $\Gamma.\;$It will be called the plactic monoid of type $G_{2}$
and denoted $Pl(G_{2})$.\ By using this monoid, we describe the corresponding
column insertion algorithm which yields a Schensted type correspondence in
$\Gamma$. Note that such a correspondence also exists for types $C$ \cite{ba}
\cite{L0}, $B$ and $D$ \cite{L2(RS)}.

The second part of this paper is devoted to the computation of the global
basis of $V(\lambda).$ We introduced a $q$-analogue $W(\Lambda_{2})$ of the
$2$-th wedge product of $V(\Lambda_{1})$.$\;$The representation $W(\Lambda
_{2})$ is not irreducible but contains an irreducible component isomorphic to
$V(\Lambda_{2}).\;$To make the notation homogeneous set $W(\Lambda
_{1})=V(\Lambda_{1}).\;$For $p=1,2,$ $W(\Lambda_{p})$ has a simple crystal
basis naturally indexed in terms of column shaped Young tableaux of height
$p.$ Then we give explicit formulas for the expansion of the global basis of
$V(\Lambda_{p}),$ $p=1,2$ on these bases.\ In the general case, we embed
$V(\lambda)$ with $\lambda=\lambda_{1}\Lambda_{1}+\lambda_{2}\Lambda_{2}$ in
\[
W(\lambda)=W(\Lambda_{1})^{\otimes\lambda_{1}}\otimes W(\Lambda_{2}%
)^{\otimes\lambda_{2}}.
\]
The tensor product of the crystal bases of the modules $W(\Lambda_{p}),$
$p=1,2$ occurring in $W(\lambda)$ is a natural basis $\{v_{\tau}\}$ of
$W(\lambda)$ indexed by combinatorial objects $\tau$ called tabloids. Then we
describe an algorithm similar to those given in \cite{L-T}, \cite{L1} and
\cite{lec3} which provides the expansion of the global basis of $V(\lambda)$
on the basis $\{v_{\tau}\}.$ Note that the coefficients of this expansion are
integral that is belong to $\mathbb{Z}[q,q^{-1}].$

\section{Background on $U_{q}(G_{2})$}

In this section we briefly review the basic facts that we shall need
concerning the representation theory of $U_{q}(G_{2})$ and the notions of
crystal basis and canonical basis of their representations. The reader is
referred to \cite{Ch-Pr}, \cite{Jan}, \cite{Ka1} and \cite{Ka2} for more details.

\subsection{The quantum enveloping algebras $U_{q}(G_{2})$}

The Dynkin diagram of the Lie algebra of type $G_{2}$ is%
\[
\overset{1}{\circ}\Lleftarrow\overset{2}{\circ}.
\]
Given a fixed indeterminate $q$ set%
\[
q_{i}=\left\{
\begin{tabular}
[c]{l}%
$q$ if $i=1$\\
$q^{3}$ if $i=2$%
\end{tabular}
\right.  \text{, }[m]_{i}=\frac{q_{i}^{m}-q_{i}^{-m}}{q_{i}-q_{i}^{-1}}\text{
and }[m]_{i}!=[m]_{i}[m-1]_{i}\cdot\cdot\cdot\lbrack1]_{i}.
\]
The quantized enveloping algebras $U_{q}(G_{2})$ is the associative algebra
over $\mathbb{Q}(q)$ generated by $e_{i},f_{i},t_{i},t_{i}^{-1},$ $i=1,2$,
subject to relations determined by the Cartan matrix $\left(
\begin{array}
[c]{cc}%
2 & -1\\
-3 & 2
\end{array}
\right)  $ of type $G_{2}$.\ Given $i\in\{1,2\}$ and $m\in\mathbb{N}$ we set
$e_{i}^{(m)}=e_{i}^{m}/[m]_{i}!$ and $f_{i}^{(m)}=f_{i}^{m}/[m]_{i}!$ For any
$i=1,2$ the subalgebra of $U_{q}(G_{2})$ generated by $e_{i},f_{i}$ and
$t_{i}$ is isomorphic to $U_{q}(\frak{sl}_{2})$ the quantum enveloping algebra
associated to $\frak{sl}_{2}.$

The representation theory of $U_{q}(G_{2})$ is closely parallel to that of
$G_{2}$. The weight lattice $P$ of $U_{q}(G_{2})$ is the $\mathbb{Z}$-lattice
generated by the fundamentals weights $\Lambda_{1},\Lambda_{2}.$ Write $P_{+}$
for the set of dominant weights of $U_{q}(G_{2})$.\ We denote by $V(\lambda)$
the irreducible finite dimensional $U_{q}(G_{2})$-module with highest weight
$\lambda\in P^{+}.$

Given two $U_{q}(G_{2})$-modules $M$ and $N$, we can define a structure of
$U_{q}(G_{2})$-module on $M\otimes N$ by putting:%
\begin{gather}
t_{i}(u\otimes v)=t_{i}u\otimes t_{i}v,\label{tensor1}\\
e_{i}(u\otimes v)=e_{i}u\otimes t_{i}^{-1}v+u\otimes e_{i}v,\label{tensor2}\\
f_{i}(u\otimes v)=f_{i}u\otimes v+t_{i}u\otimes f_{i}v. \label{tensor3}%
\end{gather}

In the sequel we need the following general lemma (see \cite{Jan} p.32).\ Let
$V(l)$ be the irreducible $U_{q}(\frak{sl}_{2})$-module of dimension $l+1$.

\begin{lemma}
\label{lem_sl2}Consider $v_{r}\in V(r)$ and $v_{s}\in V(s)$.\ Set
$t(v_{r})=q^{a}v_{r}$ with $a\in\mathbb{Z}$. Then for any integer $r$ one has:%
\[
f^{(m)}(v_{r}\otimes v_{s})=\overset{m}{\underset{k=0}{\sum}}q^{(m-k)(a-k)}%
\ f^{(k)}(v_{r})\otimes f^{(m-k)}(v_{s}).
\]
\end{lemma}

\subsection{Crystal basis and crystal graph of $U_{q}(G_{2})$}

The reader is referred to \cite{Jan} and \cite{Ka2} for basic definitions on
crystal bases and crystal graphs. Given $(L,B)$ and $(L^{\prime},B^{\prime})$
two crystal bases of the finite-dimensional $U_{q}(G_{2})$-modules $M$ and
$M^{\prime}$, $(L\otimes L^{\prime},$ $B\otimes B^{\prime})$ with $B\otimes
B^{\prime}=\{b\otimes b^{\prime};$ $b\in B,b^{\prime}\in B^{\prime}\}$ is a
crystal basis of $M\otimes M^{\prime}$. The action of $\widetilde{e}_{i}$ and
$\widetilde{f}_{i}$ on $B\otimes B^{\prime}$ is given by:%

\begin{align}
\widetilde{f_{i}}(u\otimes v)  &  =\left\{
\begin{tabular}
[c]{c}%
$\widetilde{f}_{i}(u)\otimes v$ if $\varphi_{i}(u)>\varepsilon_{i}(v)$\\
$u\otimes\widetilde{f}_{i}(v)$ if $\varphi_{i}(u)\leq\varepsilon_{i}(v)$%
\end{tabular}
\right. \label{TENS1}\\
&  \text{and}\nonumber\\
\widetilde{e_{i}}(u\otimes v)  &  =\left\{
\begin{tabular}
[c]{c}%
$u\otimes\widetilde{e_{i}}(v)$ if $\varphi_{i}(u)<\varepsilon_{i}(v)$\\
$\widetilde{e_{i}}(u)\otimes v$ if$\varphi_{i}(u)\geq\varepsilon_{i}(v)$%
\end{tabular}
\right.  \label{TENS2}%
\end{align}
where $\varepsilon_{i}(u)=\max\{k;\widetilde{e}_{i}^{k}(u)\neq0\}$ and
$\varphi_{i}(u)=\max\{k;\widetilde{f}_{i}^{k}(u)\neq0\}$. We say that a vertex
$b$ is a highest weight vertex if $\widetilde{e_{i}}(b)=0$ for $i=1,2.$ Then
$u\otimes v$ is a highest weight vertex if and only if
\begin{equation}
u\text{ is a highest weight vertex and }\varepsilon_{i}(v)\leq\varphi
_{i}(u)\text{ for }i=1,2. \label{higest_weight_vertex}%
\end{equation}

\subsection{Combinatorics of crystal graphs \label{sub_sec_combi}}

In this paragraph we recall Kang-Misra's combinatorial description of
$B(\lambda)$ \cite{KM}. It is based on the notion of tableau of type $G_{2}$
analogous to Young tableau for type $A$.

The crystal graphs of the representation $B(\Lambda_{1})$ is%
\begin{equation}
1\overset{1}{\rightarrow}2\overset{2}{\rightarrow}3\overset{1}{\rightarrow
}0\overset{1}{\rightarrow}\overline{3}\overset{2}{\rightarrow}\overline
{2}\overset{1}{\rightarrow}\overline{1}.\label{vect_G}%
\end{equation}
Set $\Gamma=\underset{l\geq0}{\bigoplus}B(\Lambda_{1})^{\otimes l}$.\ Each
vertex $x_{1}\otimes x_{2}\cdot\cdot\cdot\otimes x_{l}$ of $\Gamma$ may be
identified with the word $x_{1}\cdot\cdot\cdot x_{l}$ on the totally ordered
alphabet
\[
\mathcal{G}=\{1\prec2\prec3\prec0\prec\overline{3}\prec\overline{2}%
\prec\overline{1}\}.
\]
For any word $w$ we denote by $B(w)$ the connected component of $\Gamma$
containing $w$ and we write $l(w)$ for the length of $w.$ For $i=1,2,3$ let
$d_{i}$ be the number of letters $i$ in $w$ minus the number of letters
$\overline{i}.$ Then the weight of the vertex $w$ is given by%
\[
\mathrm{wt}(w)=(d_{1}-d_{2}+2d_{3})\Lambda_{1}+(d_{2}-d_{3})\Lambda_{2}.
\]

A column $C$ on the alphabet $\mathcal{G}$ is a Young diagram of shape column
and height $1$ or $2$ filled by letters of $\mathcal{G}$ and such that%
\[
C=%
\begin{tabular}
[c]{|l|}\hline
$a$\\\hline
\end{tabular}
\text{ with }a\in\mathcal{G}\text{ or }C=%
\begin{tabular}
[c]{|l|}\hline
$a$\\\hline
$b$\\\hline
\end{tabular}
\text{ with }a\prec b\in\mathcal{G}\text{ or }C=%
\begin{tabular}
[c]{|l|}\hline
$0$\\\hline
$0$\\\hline
\end{tabular}
.
\]
Write $h(C)$ for the height of $C.\;$Let $\mathbf{C}(p)$ , $p=1,2$ be the set
of columns of height $p.$ To each column $C,$ we associate its reading
$\mathrm{w}(C)$ that is the word obtained by reading the letters of $C$ from
top to bottom.

\noindent Denote by $\mathrm{dist}(a,b)$ the number of arrows between the
vertices $a$ and $b$ in the crystal (\ref{vect_G}) of $B(\Lambda_{1}).$ A
column $C$ is said admissible if $h(C)=1$ or $C=%
\begin{tabular}
[c]{|l|}\hline
$a$\\\hline
$b$\\\hline
\end{tabular}
$ with
\[
\left\{
\begin{tabular}
[c]{l}%
$\mathrm{dist}(a,b)\leq2$ if $a=1$ or $0$\\
$\mathrm{dist}(a,b)\leq3$ otherwise
\end{tabular}
\right.  .
\]
Given two admissible columns $C_{1}$ and $C_{2}$ we write $C_{1}\preceq C_{2}$
when the juxtaposition $C_{1}C_{2}$ of these columns verifies one of the
following assertions%
\begin{equation}
\left\{
\begin{tabular}
[c]{l}%
$\text{$\mathrm{(i)}$ }C_{1}C_{2}=$%
\begin{tabular}
[c]{|l|l|}\hline
$a$ & $b$\\\hline
\end{tabular}
with $a\preceq b\text{ and }(a,b)\neq(0,0)$\\
\ \ \ \\
$\text{$\mathrm{(ii)}$ }C_{1}C_{2}=%
\begin{tabular}
[c]{|l|l}\hline
$a$ & \multicolumn{1}{|l|}{$c$}\\\hline
$b$ & \\\cline{1-1}%
\end{tabular}
\text{ with }a\preceq c\text{ and }(a,c)\neq(0,0)$\\
\ \ \ \\
$\text{$\mathrm{(iii)}$ }C_{1}C_{2}=%
\begin{tabular}
[c]{|l|l|}\hline
$a$ & $c$\\\hline
$b$ & $d$\\\hline
\end{tabular}
$ with $\left\{
\begin{tabular}
[c]{l}%
$a\preceq c\text{ and }(a,c)\neq(0,0)$\\
$b\preceq d\text{ and }(b,d)\neq(0,0)$%
\end{tabular}
\right.  $ and $\left\{
\begin{tabular}
[c]{l}%
$\mathrm{dist}(a,d)\geq3$ if $a=2,3,0$\\
$\mathrm{dist}(a,d)\geq2$ if $a=\overline{3}$%
\end{tabular}
\right.  $%
\end{tabular}
\right.  .\label{cond_tab}%
\end{equation}

Consider $\lambda=\lambda_{1}\Lambda_{1}+\lambda_{2}\Lambda_{2}$ a dominant
weight.\ Then the vertex
\[
b_{\lambda}=(1\otimes2)^{\otimes\lambda_{2}}\otimes(1)^{\otimes\lambda_{1}}%
\]
is of highest weight $\lambda$ in $\Gamma.\;$In the sequel we will identify
$B(\lambda)$ with $B(b_{\lambda}).$ We associate to $\lambda$ the Young
diagram $Y(\lambda)$ having $\lambda_{i}$ $i=1,2$ columns of height $i$. By
definition, a tableau $T$ of type $G_{2}$ and shape $\lambda$ is a filling of
$Y(\lambda)$ by letters of $\mathcal{G}$ such that the columns $C_{i}$ of
$T=C_{1}\cdot\cdot\cdot C_{s}$ are admissible and $C_{i}\preceq C_{i+1}$ for
$i=1,...,s-1.$ We write $\mathbf{T}_{G}(\lambda)$ for the set of tableaux of
type $G_{2}$ and shape $\lambda.\;$The reading of the tableau $T=C_{1}%
\cdot\cdot\cdot C_{s}$ is the word $\mathrm{w}(T)=\mathrm{w}(C_{s})\cdot
\cdot\cdot\mathrm{w}(C_{1})$ obtained by reading the columns of $T$ from right
to left.

\begin{theorem}
(Kang-Misra) \label{KM}

\textrm{(i)}: The vertices of $B(\Lambda_{p})$ $p=1,2$ are the readings of the
admissible columns of height $p.$

\textrm{(ii)}: The vertices of $B(\lambda)$ are the readings of the tableaux
of type $G_{2}$ and shape $\lambda$.
\end{theorem}%

\begin{gather}%
\begin{tabular}
[c]{lllllllllllll}%
11 & $\overset{1}{\rightarrow}$ & 21 & $\overset{2}{\rightarrow}$ & 31 &
$\overset{1}{\rightarrow}$ & 01 & $\overset{1}{\rightarrow}$ & \={3}1 &
$\overset{2}{\rightarrow}$ & \={2}1 & $\overset{1}{\rightarrow}$ & \={1}1\\
&  & {\scriptsize 1}$\downarrow$ &  &  &  &  &  & {\scriptsize 1}$\downarrow$%
&  &  &  & {\scriptsize 1}$\downarrow$\\
12 &  & 22 & $\overset{2}{\rightarrow}$ & 32 & $\overset{1}{\rightarrow}$ &
02 &  & \={3}2 & $\overset{2}{\rightarrow}$ & \={2}2 &  & \={1}2\\
{\scriptsize 2}$\downarrow$ &  &  &  & {\scriptsize 2}$\downarrow$ &  &
{\scriptsize 2}$\downarrow$ &  &  &  & {\scriptsize 2}$\downarrow$ &  &
{\scriptsize 2}$\downarrow$\\
13 & $\overset{1}{\rightarrow}$ & 23 &  & 33 & $\overset{1}{\rightarrow}$ &
03 & $\overset{1}{\rightarrow}$ & \={3}3 &  & \={2}3 & $\overset
{1}{\rightarrow}$ & \={1}3\\
&  & {\scriptsize 1}$\downarrow$ &  &  &  &  &  & {\scriptsize 1}$\downarrow$%
&  &  &  & {\scriptsize 1}$\downarrow$\\
10 &  & 20 & $\overset{2}{\rightarrow}$ & 30 & $\overset{1}{\rightarrow}$ &
00 &  & 30 & $\overset{2}{\rightarrow}$ & \={2}0 &  & \={1}\={3}\\
{\scriptsize 1}$\downarrow$ &  & {\scriptsize 1}$\downarrow$ &  &  &  &
{\scriptsize 1}$\downarrow$ &  & {\scriptsize 1}$\downarrow$ &  &
{\scriptsize 1}$\downarrow$ &  & {\scriptsize 1}$\downarrow$\\
1\={3} &  & 2\={3} & $\overset{2}{\rightarrow}$ & 3\={3} &  & 0\={3} &  &
\={3}\={3} & $\overset{2}{\rightarrow}$ & \={2}\={3} &  & \={1}\={3}\\
{\scriptsize 2}$\downarrow$ &  &  &  & {\scriptsize 2}$\downarrow$ &  &
{\scriptsize 2}$\downarrow$ &  &  &  & {\scriptsize 2}$\downarrow$ &  &
{\scriptsize 2}$\downarrow$\\
1\={2} & $\overset{1}{\rightarrow}$ & 2\={2} &  & 3\={2} & $\overset
{1}{\rightarrow}$ & 0\={2} & $\overset{1}{\rightarrow}$ & \={3}\={2} &  &
\={2}\={2} & $\overset{1}{\rightarrow}$ & \={1}\={2}\\
&  & {\scriptsize 1}$\downarrow$ &  &  &  &  &  & {\scriptsize 1}$\downarrow$%
&  &  &  & {\scriptsize 1}$\downarrow$\\
1\={1} &  & 2\={1} & $\overset{2}{\rightarrow}$ & 3\={1} & $\overset
{1}{\rightarrow}$ & 0\={1} &  & \={3}\={1} & $\overset{2}{\rightarrow}$ &
\={2}\={1} &  & \={1}\={1}%
\end{tabular}
\label{carre_tens}\\
\text{the crystal }B(\Lambda_{1})^{\otimes2}\simeq B(2\Lambda_{1})%
%TCIMACRO{\tbigoplus }%
%BeginExpansion
{\textstyle\bigoplus}
%EndExpansion
B(\Lambda_{2})%
%TCIMACRO{\tbigoplus }%
%BeginExpansion
{\textstyle\bigoplus}
%EndExpansion
B(\Lambda_{1})%
%TCIMACRO{\tbigoplus }%
%BeginExpansion
{\textstyle\bigoplus}
%EndExpansion
B(0)\nonumber
\end{gather}

\subsection{Canonical basis of $V(\lambda)$}

Denote by $F\mapsto\overline{F}$ the involution of $U_{q}(G_{2})$ defined as
the ring automorphism satisfying%
\[
\overline{q}=q^{-1},\text{ \ }t_{i}=t_{i}^{-1},\text{ \ \ }\overline{e_{i}%
}=e_{i},\text{ \ \ }\overline{f_{i}}=f_{i}\text{ \ \ \ for }i=1,2.
\]
By writing each vector $v$ of $V(\lambda)$ in the form $v=Fv_{\lambda}$ where
$F\in U_{q}(G_{2})$, we obtain an involution of $V(\lambda)$ defined by%
\[
\overline{v}=\overline{F}v_{\lambda}.
\]
Let $U_{\mathbb{Q}}^{-}$ be the subalgebra of $U_{q}(G_{2})$ generated over
$\mathbb{Q}[q,q^{-1}]$ by the $f_{i}^{(k)}$ and set $V_{\mathbb{Q}}%
(\lambda)=U_{\mathbb{Q}}^{-}v_{\lambda}$. We can now state:

\begin{theorem}
\label{Th_K}(Kashiwara) There exists a unique $\mathbb{Q[}q,q^{-1}]$-basis
$\{G(T);$ $T\in\mathbf{T}_{G}(\lambda)\}$ of $V_{\mathbb{Q}}(\lambda)$ such
that:%
\begin{gather}
G(T)\equiv\mathrm{w(}T)\text{ }\mathrm{mod}\text{ }qL(\lambda
),\label{cond_cong}\\
\overline{G(T)}=G(T).\label{cond_invo}%
\end{gather}
\end{theorem}

\noindent Note that $G(T)\in V_{\mathbb{Q}}(\lambda)\cap L(\lambda)$.\ The
basis $\{G(T);T\in\mathbf{T}_{G}(\lambda)\}$ is called the lower global (or
canonical) basis of $V(\lambda)$

\section{Plactic monoid and Schensted's type correspondence}

\subsection{The plactic monoid}

\begin{definition}
\label{def_sam_plac}Let $w_{1}$ and $w_{2}$ be two words on $\mathcal{G}$ . We
write $w_{1}\sim w_{2}$ when these two words occur at the same place in two
isomorphic connected components of the crystal $\Gamma$.
\end{definition}

From Theorem \ref{KM} it follows that for any word $w\in\mathcal{G}^{\ast}$
there exists a unique tableau of type $G_{2}$ $P(w)$ such that $w\sim
\mathrm{w(}P(w))$. So the set $\mathcal{G}^{\ast}/\sim$ can be identified with
the set of tableaux of type $G_{2}$. Our aim is now to show that $\sim$ are in
fact a congruences $\equiv$ so that $\mathcal{G}^{\ast}/\sim$ is in a natural
way endowed with a multiplication.

\noindent Set
\[
S=\{21,31,01,\overline{3}1,\overline{3}2,\overline{2}1,\overline{2}%
2,\overline{1}1,\overline{1}2,\overline{2}3,\overline{1}3,\overline
{1}0,\overline{1}\ \overline{3},\overline{1}\ \overline{2}\}.
\]
To defined the plactic relations for type $G_{2}$ we need the bijection
$\Theta$ from $S$ to $B(12)$ defined by%
\begin{equation}%
\begin{tabular}
[c]{|c|c|c|c|c|c|c|c|c|c|c|c|c|c|c|}%
$w$ & $%
\begin{tabular}
[c]{l}
\ \ \\
\ \
\end{tabular}
21$ & $31$ & $01$ & $\bar{3}1$ & $\bar{3}2$ & $\bar{2}1$ & $\bar{2}2$ &
$\bar{1}1$ & $\bar{1}2$ & $\bar{2}3$ & $\bar{1}3$ & $\bar{1}0$ & $\bar{1}%
\bar{3}$ & $\bar{1}\bar{2}$\\\hline
$\Theta(w)$ &
\begin{tabular}
[c]{l}
\ \ \\
\ \
\end{tabular}
$12$ & $13$ & $23$ & $20$ & $2\bar{3}$ & $30$ & $3\bar{3}$ & $00$ & $0\bar{3}$%
& $3\bar{2}$ & $0\bar{2}$ & $\bar{3}\bar{2}$ & $\bar{3}\bar{1}$ & $\bar{2}%
\bar{1}$%
\end{tabular}
\label{tab-cores}%
\end{equation}

\begin{definition}
\label{def_monoi}The monoid $Pl(G_{2})$ is the quotient of the free monoid
$\mathcal{G}^{\ast}$ by the relations:%
\begin{equation}
10\equiv1,\text{ }1\overline{3}\equiv2,\text{ }1\overline{2}\equiv3,\text{
}2\overline{2}\equiv0,\text{ }0\overline{1}\equiv\overline{1},\text{
}3\overline{1}\equiv\overline{2},\text{ }2\overline{1}\equiv\overline
{3}.\tag{$R_1$}\label{R1}%
\end{equation}%

\begin{equation}
1\overline{1}\equiv\emptyset.\tag{$R_2$}\label{R2}%
\end{equation}%

\begin{equation}
abc\equiv\left\{
\begin{tabular}
[c]{l}%
$a\Theta(bc)$ if $bc\in S$\\
$\Theta^{-1}(ab)c$ otherwise
\end{tabular}
\right.  \text{ with }ab\in B(12)\text{ and }bc\in B(11).\tag{$R_3$}\label{R3}%
\end{equation}%

\begin{equation}
xyz\equiv\Theta^{-1}(xy)z\text{ with }xy\in B(\Lambda_{2})\text{ and }yz\in
B(\Lambda_{2}).\tag{$R_4$}\label{R4}%
\end{equation}
\end{definition}

For any word $w$ occurring in the left hand side of a relation $R_{i},$
$i=1,...,4$ we write $\xi_{i}(w)$ the word occurring in the right hand side of
this relation.

\begin{proposition}
\label{prop_iso_plactic}The maps $\xi_{1},\xi_{2},\xi_{3}$ and $\xi_{4}$ are
respectively the crystal graph isomorphisms
\[
\mathrm{(i):}\text{ }B(10)\overset{\sim}{\rightarrow}B(1),\text{
}\mathrm{(ii):}\text{ }B(1\overline{1})\overset{\sim}{\rightarrow}%
B(\emptyset),\text{ }\mathrm{(iii):}\text{ }B(121)\overset{\sim}{\rightarrow
}B(112)\text{ and }\mathrm{(iv)}\text{: }B(123)\overset{\sim}{\rightarrow
}B(110).
\]
\end{proposition}

\begin{proof}
The simplest way to prove this proposition consists in the computation of the
crystals above.\ The crystals $B(10)$ and $B(1\overline{1})$ occur in
(\ref{carre_tens}) and we see that $\xi_{1}$ and $\xi_{2}$ are the crystal
graph isomorphisms \textrm{(i) }and \textrm{(ii).\ }We obtain that $\xi_{4}$
is the crystal graph isomorphism \textrm{(iv) }from the crystals
(\ref{B(123)}) and (\ref{B(110)}) below.\ Similarly we can prove that $\xi
_{3}$ is the crystal graph isomorphism between $B(121)$ and $B(112)$ by
drawing these two crystals which contain $64$ vertices.
\end{proof}

\noindent\textbf{Remark }Write $\xi_{1}^{\prime}$ for the crystal graph
isomorphism $B(110)\overset{\sim}{\rightarrow}B(11)$, then $\xi_{4}^{\prime
}=\xi_{1}^{\prime}\xi_{4}$ is the crystal isomorphism $B(123)\overset{\sim
}{\rightarrow}B(11).$ By replacing Relation $R_{4}$ in the definition of
$Pl(G_{2})$ above by a relation $R_{4}^{\prime}$ describing the isomorphism
$\xi_{4}^{\prime}$ we obtain Littelmann's presentation of $Pl(G_{2})$
\cite{lit}.\ The two presentations coincide excepted for the description of
$\xi_{3}$ which is not totally exact in \cite{lit}.%

\begin{align}
&  {%
\begin{tabular}
[c]{lllllllllllll}%
123 & $\overset{1}{\rightarrow}$ & 120 & $\overset{2}{\rightarrow}$ & 130 &
$\overset{1}{\rightarrow}$ & 230 & $\overset{1}{\rightarrow}$ & 200 &
$\overset{2}{\rightarrow}$ & 300 & $\overset{1}{\rightarrow}$ & 000\\
&  & {\scriptsize 1}$\downarrow$ &  &  &  &  &  & {\scriptsize 1}$\downarrow$%
&  &  &  & {\scriptsize 1}$\downarrow$\\
&  & 12\={3} & $\overset{2}{\rightarrow}$ & 13\={3} & $\overset{1}%
{\rightarrow}$ & 23\={3} &  & 20\={3} & $\overset{2}{\rightarrow}$ & 30\={3} &
& 00\={3}\\
&  &  &  & {\scriptsize 2}$\downarrow$ &  & {\scriptsize 2}$\downarrow$ &  &
&  & {\scriptsize 2}$\downarrow$ &  & {\scriptsize 2}$\downarrow$\\
&  &  &  & 13\={2} & $\overset{1}{\rightarrow}$ & 23\={2} & $\overset
{1}{\rightarrow}$ & 20\={2} &  & 30\={2} & $\overset{1}{\rightarrow}$ &
00\={2}\\
&  &  &  &  &  &  &  & {\scriptsize 1}$\downarrow$ &  &  &  & {\scriptsize 1}%
$\downarrow$\\
&  &  &  &  &  &  &  & 2\={3}\={2} & $\overset{2}{\rightarrow}$ & 3\={3}\={2}%
&  & 0\={3}\={2}\\
&  &  &  &  &  &  &  & {\scriptsize 1}$\downarrow$ &  & {\scriptsize 1}%
$\downarrow$ &  & {\scriptsize 1}$\downarrow$\\
&  &  &  &  &  &  &  & 2\={3}\={1} & $\overset{2}{\rightarrow}$ & 3\={3}\={1}%
&  & 0\={3}\={1}\\
&  &  &  &  &  &  &  &  &  & {\scriptsize 2}$\downarrow$ &  & {\scriptsize 2}%
$\downarrow$\\
&  &  &  &  &  &  &  &  &  & 3\={2}\={1} & $\overset{1}{\rightarrow}$ &
0\={2}\={1}\\
&  &  &  &  &  &  &  &  &  &  &  & {\scriptsize 1}$\downarrow$\\
&  &  &  &  &  &  &  &  &  &  &  & \={3}\={2}1
\end{tabular}
}\label{B(123)}\\
&  {%
\begin{tabular}
[c]{lllllllllllll}%
110 & $\overset{1}{\rightarrow}$ & 210 & $\overset{2}{\rightarrow}$ & 310 &
$\overset{1}{\rightarrow}$ & 010 & $\overset{1}{\rightarrow}$ & \={3}10 &
$\overset{2}{\rightarrow}$ & \={2}10 & $\overset{1}{\rightarrow}$ & \={1}10\\
&  & {\scriptsize 1}$\downarrow$ &  &  &  &  &  & {\scriptsize 1}$\downarrow$%
&  &  &  & {\scriptsize 1}$\downarrow$\\
&  & 21\={3} & $\overset{2}{\rightarrow}$ & 31\={3} & $\overset{1}%
{\rightarrow}$ & 01\={3} &  & \={3}1\={3} & $\overset{2}{\rightarrow}$ &
\={2}1\={3} &  & \={1}1\={3}\\
&  &  &  & {\scriptsize 2}$\downarrow$ &  & {\scriptsize 2}$\downarrow$ &  &
&  & {\scriptsize 2}$\downarrow$ &  & {\scriptsize 2}$\downarrow$\\
&  &  &  & 31\={2} & $\overset{1}{\rightarrow}$ & 01\={2} & $\overset
{1}{\rightarrow}$ & \={3}1\={2} &  & \={2}1\={2} & $\overset{1}{\rightarrow}$%
& \={1}1\={2}\\
&  &  &  &  &  &  &  & {\scriptsize 1}$\downarrow$ &  &  &  & {\scriptsize 1}%
$\downarrow$\\
&  &  &  &  &  &  &  & \={3}2\={2} & $\overset{2}{\rightarrow}$ & \={2}2\={2}%
&  & 12\={2}\\
&  &  &  &  &  &  &  & {\scriptsize 1}$\downarrow$ &  & {\scriptsize 1}%
$\downarrow$ &  & {\scriptsize 1}$\downarrow$\\
&  &  &  &  &  &  &  & \={3}2\={1} & $\overset{2}{\rightarrow}$ & \={2}2\={1}%
&  & 12\={1}\\
&  &  &  &  &  &  &  &  &  & {\scriptsize 2}$\downarrow$ &  & {\scriptsize 2}%
$\downarrow$\\
&  &  &  &  &  &  &  &  &  & \={2}3\={1} & $\overset{1}{\rightarrow}$ &
13\={1}\\
&  &  &  &  &  &  &  &  &  &  &  & {\scriptsize 1}$\downarrow$\\
&  &  &  &  &  &  &  &  &  &  &  & 0\={1}1
\end{tabular}
} \label{B(110)}%
\end{align}

\begin{theorem}
\label{th_good_rela}Given two words $w_{1}$ and $w_{2}$%
\begin{equation}
w_{1}\sim w_{2}\Longleftrightarrow w_{1}\equiv w_{2}\label{good_rela}%
\end{equation}
\end{theorem}

\noindent To prove this, we need the following lemma.

\begin{lemma}
\label{lem_aux}Let $w$ be a highest weight vertex of $\Gamma$.\ Then
$\mathrm{w(}P(w))\equiv w$.
\end{lemma}

\begin{proof}
We proceed by induction on $l(w)$. If $l(w)=1$ then $w=1$ and $P(w)=%
\begin{tabular}
[c]{|l|}\hline
$1$\\\hline
\end{tabular}
$. Suppose the proposition true for the highest weight vertices of length
$\leq l$ and consider a highest weight vertex $w$ of length $l+1$. Write
$w=vx$ where $x\in\mathcal{G}$ and $v$ is a word of length $l$. Then if
follows from (\ref{higest_weight_vertex}) that $v$ is a highest weight vertex
and $\varepsilon_{i}(x)\leq\varphi_{i}(v)$ for $i=1,2.\;$So by induction
$v\equiv\mathrm{w(}P(v))$. Hence $w\equiv\mathrm{w(}P(v))x$. Let $P^{\prime}$
be the tableau obtained by erasing the first column of $P(v).\;$The condition
$\varepsilon_{i}(x)\leq\varphi_{i}(v)$ implies that only the following
situations can occur:

\begin{enumerate}
\item $x=1$ and $P(w)$ is obtained by adding a column of height $1$ in $P(v)$.
We have $w\equiv\mathrm{w(}P(v))1\equiv\mathrm{w(}P(w))$ Indeed, by relation
$R_{3},$ $1$ commutes in $Pl(G_{2})$ with all the column words $12$ occurring
in $\mathrm{w(}P(v)).$

\item $x=2$ and $P(w)$ is obtained from $P(v)$ by erasing a column of height
$1$ and adding a column of height $2$. Similarly we have $w\equiv
\mathrm{w(}P(w))$ because in $Pl(G_{2})$, $2$ commute with all the columns
words $12$ occurring in $\mathrm{w(}P(v)).$

\item $x=3,$ $P(v)$ has at least a column of height $2$ and $P(w)$ is obtained
from $P(v)$ by erasing a column of height $2$ and adding two columns of height
$1$. We can write $\mathrm{w}(P(v))=\mathrm{w(}P^{\prime})12$.\ Then
$w\equiv\mathrm{w(}P^{\prime})123\equiv\mathrm{w(}P^{\prime})11$ by applying
relations $R_{4}$ and $R_{1}$.$\;$We obtain $w\equiv\mathrm{w(}P(w))$ by using
the same argument than in case $1.$

\item $x=0,$ $P(v)$ has at least a column of height $1$ and $P(w)=P(v).$ If
$\mathrm{w}(P(v))=\mathrm{w(}P^{\prime})1$ then $w\equiv\mathrm{w(}P^{\prime
})10\equiv\mathrm{w(}P^{\prime})1$ by relation $R_{1}.$ If $\mathrm{w}%
(P(v))=\mathrm{w(}P^{\prime})12$ then $w\equiv\mathrm{w(}P^{\prime}%
)120\equiv\mathrm{w(}P^{\prime})21$ by relations $R_{1}$ and $R_{4}.$ In the
both case $w\equiv\mathrm{w(}P(w)).$

\item $x=\overline{3}$, $P(v)$ has at least two columns of height $1$ and
$P(w)$ is obtained from $P(v)$ by erasing two columns of height $1$ and adding
a column of height $2.$ If $\mathrm{w}(P(v))=\mathrm{w(}P^{\prime})1$ then
$w\equiv\mathrm{w(}P^{\prime})1\overline{3}\equiv\mathrm{w(}P^{\prime})2$ by
applying relation $R_{1}.$ If $\mathrm{w}(P(v))=\mathrm{w(}P^{\prime})12$ then
$w\equiv\mathrm{w(}P^{\prime})12\overline{3}\equiv\mathrm{w(}P^{\prime})22$ by
applying relations $R_{4}$ and $R_{1}.$ In the both cases $w\equiv\mathrm{w}(P(w)).$

\item $x=\overline{2}$, $P(v)$ has at least a column of height $2$ and $P(w)$
is obtained from $P(v)$ by erasing a column of height $2$ and adding a column
of height $1.$ We can write $\mathrm{w}(P(v))=\mathrm{w(}P^{\prime})12.\;$So
$w\equiv\mathrm{w(}P^{\prime})12\overline{2}\equiv\mathrm{w(}P^{\prime
})1\equiv\mathrm{w(}P(w))$ by applying relation $R_{1}$ two times.

\item $x=\overline{1},$ $P(v)$ has at least a column of height $1$ and $P(w)$
is obtained from $P(v)$ by erasing a column of height $1.\;$If $\mathrm{w}%
(P(v))=\mathrm{w(}P^{\prime})1$ then $w\equiv\mathrm{w(}P^{\prime}%
)1\overline{1}\equiv\mathrm{w(}P^{\prime})$ by applying relation $R_{2}.$ If
$\mathrm{w}(P(v))=\mathrm{w(}P^{\prime})12$ then $w\equiv\mathrm{w(}P^{\prime
})12\overline{1}\equiv\mathrm{w(}P^{\prime})2$ by applying relation $R_{1}$
two times$.$ In the both cases $w\equiv\mathrm{w}(P(w)).$
\end{enumerate}
\end{proof}

\begin{proof}
(of Theorem \ref{th_good_rela}) By Proposition \ref{prop_iso_plactic} and
(\ref{TENS1}), (\ref{TENS2}), the plactic relations of Definition
\ref{def_monoi} are compatible with Kashiwara's operators, that is, for any
words $w_{1}$ and $w_{2}$ such that $w_{1}\equiv w_{2}$ one has:%
\begin{equation}
\left\{
\begin{tabular}
[c]{l}%
$\widetilde{e}_{i}(w_{1})\equiv\widetilde{e}_{i}(w_{2})\text{ and }%
\varepsilon_{i}(w_{1})=\varepsilon_{i}(w_{2})$\\
$\widetilde{f}_{i}(w_{1})\equiv\widetilde{f}_{i}(w_{2})\text{ and }\varphi
_{i}(w_{1})=\varphi_{i}(w_{2}).$%
\end{tabular}
\right.  \label{fonda_compatib}%
\end{equation}
Hence:%
\[
w_{1}\equiv w_{2}\Longrightarrow w_{1}\sim w_{2}.
\]

From Lemma \ref{lem_aux}, we obtain that two highest weight vertices
$w_{1}^{0}$ and $w_{2}^{0}$ with the same weight $\lambda$ verify $w_{1}%
^{0}\equiv w_{2}^{0}$. Indeed there is only one tableau of type $G_{2}$ whose
reading is a highest vertex of weight $\lambda$. Now suppose that $w_{1}\sim
w_{2}$ and denote by $w_{1}^{0}$ and $w_{2}^{0}$ the highest weight vertices
of $B(w_{1})$ and $B(w_{2})$. We have $w_{1}^{0}\equiv w_{2}^{0}$. Set
$w_{1}=\widetilde{F}\,w_{1}^{0}$ where $\widetilde{F}$ is a product of
Kashiwara's operators $\widetilde{f}_{i}$, $i=1,2$. Then $w_{2}=\widetilde
{F}\,w_{2}^{0}$ because $w_{1}\sim w_{2}$. So by (\ref{fonda_compatib}) we
obtain%
\[
w_{1}^{0}\equiv w_{2}^{0}\Longrightarrow\widetilde{F}\,w_{1}^{0}%
\equiv\widetilde{F}\,w_{2}^{0}\Longrightarrow w_{1}\equiv w_{2}.
\]
\end{proof}

\subsection{Bumping algorithm}

Now we are going to see how the orthogonal tableau $P(w)$ may be computed for
each vertex $w$ by using an insertion scheme analogous to bumping algorithm
for type $A$. As a first step, we describe $P(w)$ when $w=\mathrm{w}(C)x$,
where $x$ and $C$ are respectively a letter and an admissible column. This
will be called `` the insertion of the letter $x$ in the admissible column $C$
'' and denoted by $x\rightarrow C$. Then we will be able to obtain $P(w)$ when
$w=\mathrm{w}(T)x$ with $x$ a letter and $T$ an orthogonal tableau. This will
be called `` the insertion of the letter $x$ in the orthogonal tableau $T$ ''
and denoted by $x\rightarrow T$. Our construction of $P$ will be recursive, in
the sense that if $P(u)=T$ and $x$ is a letter, then $P(ux)=x\rightarrow T$.

\subsubsection{\label{x_in_C}Insertion of a letter in an admissible column}

\noindent When $h(C)=1$ and $C=$%
\begin{tabular}
[c]{|l|}\hline
$a$\\\hline
\end{tabular}
we have%
\begin{equation}
x\rightarrow C=\left\{
\begin{tabular}
[c]{l}%
$\mathrm{(i):}$ $%
\begin{tabular}
[c]{|l|l|}\hline
$a$ & $x$\\\hline
\end{tabular}
\ $if $ax\in B(11)$\\
\ \\
$\mathrm{(ii):}$ $%
\begin{tabular}
[c]{|l|}\hline
$a$\\\hline
$x$\\\hline
\end{tabular}
\ $if $ax\in B(12)$\\
\ \\
$\mathrm{(iii):}$ $%
\begin{tabular}
[c]{|l|}\hline
$a^{\prime}$\\\hline
\end{tabular}
\ $with $a^{\prime}=\xi_{1}(ax)$ if $ax\in B(10)$\\
\ \\
$\mathrm{(iv):}$ $\emptyset$ if $ax=1\overline{1}$%
\end{tabular}
\right.  .\label{insert_row}%
\end{equation}
Indeed in each cases $\mathrm{(i)}$ to $\mathrm{(iv)}$ $x\rightarrow C$ is a
tableau of type $G_{2}$ such that $\mathrm{w}(x\rightarrow C)\equiv\mathrm{w}(C)x.$

\noindent When $h(C)=2$ and $C=$%
\begin{tabular}
[c]{|l|}\hline
$a$\\\hline
$b$\\\hline
\end{tabular}
we have
\[
x\rightarrow C=\left\{
\begin{tabular}
[c]{l}%
$\mathrm{(v):}$
\begin{tabular}
[c]{|l|l}\hline
$a^{\prime}$ & \multicolumn{1}{|l|}{$x^{\prime}$}\\\hline
$b^{\prime}$ & \\\cline{1-1}%
\end{tabular}
with $x^{\prime}a^{\prime}b^{\prime}=\xi_{3}(abx)$ if $bx$ is not a column
word\\
\ \ \ \\
$\mathrm{(vi):}$
\begin{tabular}
[c]{|l|l|}\hline
$x^{\prime}$ & $y^{\prime}$\\\hline
\end{tabular}
with $a^{\prime}x^{\prime}=\xi_{1}\xi_{4}(abx)$ if $bx$ is an admissible
column word\\
\ \ \ \ \\
$\mathrm{(vii):}$
\begin{tabular}
[c]{|l|}\hline
$x^{\prime}$\\\hline
\end{tabular}
with $x^{\prime}=\xi_{1}(a\xi_{1}(bx))$ if $bx$ is a non admissible column
word
\end{tabular}
\right.  .
\]
Indeed in cases $\mathrm{(v)}$ and $\mathrm{(vi)}$ $x\rightarrow C$ is a
tableau of type $G_{2}$ such that $\mathrm{w}(x\rightarrow C)\equiv
\mathrm{w}(C)x$ by Proposition \ref{prop_iso_plactic}. In case $\mathrm{(vii)}%
$, we obtain by (\ref{TENS2}) that the highest weight vertex of $B(abx)$ may
be written $12u$ with $u$ a letter such that $\varepsilon_{1}(u)=0$ and
$\varepsilon_{2}(u)\leq1.$ So $u\in\{1,3,\overline{2}\}.\;$We have
$u=\overline{2}$ otherwise $B(abx)=B(121)$ and $x\preceq b,$ or
$B(abx)=B(123)$ and $bx$ is an admissible column word.\ Hence
$B(abx)=B(12\overline{2}).$ We have%
\[
B(12\overline{2}):12\overline{2}\overset{1}{\rightarrow}12\overline{1}%
\overset{2}{\rightarrow}13\overline{1}\overset{1}{\rightarrow}23\overline
{1}\overset{1}{\rightarrow}20\overline{1}\overset{2}{\rightarrow}%
30\overline{1}\overset{1}{\rightarrow}00\overline{1}%
\]
and it is easy to verify that $\xi_{1}(a\xi_{1}(bx))$ is the image of $abx$ by
the crystal isomorphism $B(12\overline{2})\overset{\sim}{\rightarrow}B(1)$.

\noindent In cases $\mathrm{(iii),}$ $\mathrm{(iv),}$ $\mathrm{(vi)}$ and
$\mathrm{(vii)}$ we have $l(x\rightarrow C)<l(\mathrm{w}(C)x)$.\ We will say
that the insertion procedure causes a contraction. Note that if the words
$\mathrm{w}(C_{1})x_{1}$ and $\mathrm{w}(C_{2})x_{2}$ (where $C_{1}$, $C_{2}$
are admissible columns and $x_{1},$ $x_{2}$ are letters) belongs to the same
connected component, the insertions $x_{1}\rightarrow C_{1}$ and
$x_{2}\rightarrow C_{2}$ are of the same type $\mathrm{(i)}$ to
$\mathrm{(vii).}$

\subsubsection{Insertion of a letter in a tableau of type $G_{2}$}

Consider a tableau $T=C_{1}C_{2}\cdot\cdot\cdot C_{r}$ of type $G_{2}$.\ The
insertion $x\rightarrow T$ is characterized by the following proposition.

\begin{proposition}
Set $T^{\prime}=C_{2}\cdot\cdot\cdot C_{s}$.

\begin{enumerate}
\item  If the insertion $x\rightarrow C_{1}$ is of type $\mathrm{(i),}$
$\mathrm{(ii)}$ or $\mathrm{(iv)}$ in \ref{x_in_C} then $x\rightarrow
T=(x\rightarrow C_{1})T^{\prime}$.

\item  If the insertion $x\rightarrow C_{1}$ is of type $\mathrm{(v)}$ in
\ref{x_in_C} with $C_{1}=$%
\begin{tabular}
[c]{|l|}\hline
$a$\\\hline
$b$\\\hline
\end{tabular}
and $C_{1}^{\prime}=$%
\begin{tabular}
[c]{|l|}\hline
$a^{\prime}$\\\hline
$b^{\prime}$\\\hline
\end{tabular}
then $x\rightarrow T=C_{1}^{\prime}(x^{\prime}\rightarrow T^{\prime})$ and is
obtained by computing successively insertions of type $\mathrm{(v)}$.

\item  If the insertion $x\rightarrow C_{1}$ is of type $\mathrm{(iii)}$ in
\ref{x_in_C} and $x\rightarrow C_{1}=$%
\begin{tabular}
[c]{|l|}\hline
$a^{\prime}$\\\hline
\end{tabular}
then $x\rightarrow T=a^{\prime}\rightarrow T^{\prime}.$

\item  If the insertion $x\rightarrow C_{1}$ is of type $\mathrm{(vi)}$ in
\ref{x_in_C}, then $x\rightarrow T=x^{\prime}\rightarrow(y^{\prime}\rightarrow
T^{\prime})$. Moreover the insertion of the letters $x^{\prime}$ and
$y^{\prime}$ in $T^{\prime}$ does not cause a new contraction.

\item  If the insertion $x\rightarrow C_{1}$ is of type $\mathrm{(vii)}$ in
\ref{x_in_C}, then $x\rightarrow T=x^{\prime}\rightarrow T^{\prime}$. Moreover
the insertion of the letter $x^{\prime}$ in $T^{\prime}$ does not cause a new contraction.
\end{enumerate}
\end{proposition}

\noindent The insertion procedure terminates because in cases $2,$ $3,4$ and
$5$ we are reduced to the insertion of a letter in a tableau whose number of
boxes is strictly less than that of $T$.

\begin{proof}
It follows from the proof of lemma \ref{lem_aux} that the difference between
the number of boxes in the shapes of $T\;$and $x\rightarrow T$ belongs to
$\{-1,0,1\}.\;$So only one contraction can occur during the insertion
$x\rightarrow T.$

\noindent In case $1$, $T$ contains only columns of height $1$.\ Hence the
proposition follows immediately from the definition of the insertions of type
$\mathrm{(i),}$ $\mathrm{(ii),}$ and $\mathrm{(iv).}$

\noindent In case $2,$ we can write by (\ref{higest_weight_vertex})
$w^{0}=v^{0}x^{0}$ where $v^{0}$ is the highest weight vertex of
$\mathrm{w}(T)$ and $x^{0}$ is a letter.\ The word $v^{0}$ is the reading of
$T^{0}$ a tableau of type $G_{2}.$ Write $T^{0}=C_{1}^{0}\cdot\cdot\cdot
C_{r}^{0}$. The insertion $x^{0}\rightarrow C_{1}^{0}$ is of type
$\mathrm{(v)}$ since it is true for the insertion $x\rightarrow C.\;$Hence
$x^{0}\in\{1,2\}$ since $\mathrm{w}(C_{1}^{0})=12.$ Let $\widehat{T}^{0}$ be
the tableau obtained from $T^{0}$ by adding a box containing the letter
$x^{0}$.\ Denote by $C_{k}^{0}$ the column of $T^{0}$ where this new box
appears.\ There exists a unique sequence $w_{0}=\mathrm{w}(T)x,...,w_{k-1}$
such that $\mathrm{w}(\widehat{T}^{0})$ is the highest weight vertex of
$B(w_{k-1})$ and such that $w_{j-1}=\mathrm{w}(C_{r})\cdot\cdot\cdot
\mathrm{w}(C_{j})x_{j-1}\mathrm{w}(C_{j-1}^{\prime})\cdot\cdot\cdot
\mathrm{w}(C_{1}^{\prime})$ is transformed into the congruent word
$w_{j}=\mathrm{w}(C_{r})\cdot\cdot\cdot\mathrm{w}(C_{j+1})x_{j}\mathrm{w}%
(C_{j}^{\prime})\cdot\cdot\cdot\mathrm{w}(C_{1}^{\prime})$ where $x_{j}$ and
$C_{j}^{\prime}$ are determined by the insertion of type $\mathrm{(v)}$
$x_{j-1}\rightarrow C_{j}=C_{j}^{\prime}%
\begin{tabular}
[c]{|l|}\hline
$x_{j}$\\\hline
\end{tabular}
$.\ The word $w_{k-1}$ is the reading of a tableau of type $G_{2}$ and is
obtained by computing $k-1$ insertions of type $\mathrm{(v)}$ so is equal to
$\mathrm{w}(x\rightarrow T).$

The rest of the proof is obtained by induction on $l(\mathrm{w}(T)).\;$If
$l(\mathrm{w}(T))=1$, the proposition follows from (\ref{insert_row}). Now
suppose the proposition true for any tableau $T$ such that $l(\mathrm{w}%
(T))\leq m$ with $m\in\mathbb{N}$.\ Then consider a tableau $T$ with
$l(\mathrm{w}(T))=m+1$ and a letter $x.$

\noindent In case $3$, $x\rightarrow T=a^{\prime}\rightarrow T^{\prime}$
since, by induction, $a^{\prime}\rightarrow T^{\prime}$ is a tableau and
$\mathrm{w}(a^{\prime}\rightarrow T^{\prime})\equiv\mathrm{w}(T)x.$

\noindent In cases $4$ and $5$ the result is obtained by using the induction
hypothesis as in case $2.$
\end{proof}

\begin{example}
Consider the tableau of type $G_{2}$%
\[
T=%
\begin{tabular}
[c]{|l|l|l|}\hline
$\mathtt{2}$ & $\mathtt{0}$ & $\mathtt{\bar{3}}$\\\hline
$\mathtt{0}$ & $\mathtt{\bar{2}}$ & $\mathtt{\bar{1}}$\\\hline
\end{tabular}
.
\]
Then the insertion $\overline{2}\rightarrow T$ may be computed as follows:%
\begin{multline*}
\overline{2}\rightarrow%
\begin{tabular}
[c]{|l|l|l|}\hline
$\mathtt{2}$ & $\mathtt{0}$ & $\mathtt{\bar{3}}$\\\hline
$\mathtt{0}$ & $\mathtt{\bar{2}}$ & $\mathtt{\bar{1}}$\\\hline
\end{tabular}
=3\rightarrow\left(  \overline{3}\rightarrow%
\begin{tabular}
[c]{|l|l|}\hline
$\mathtt{0}$ & $\mathtt{\bar{3}}$\\\hline
$\mathtt{\bar{2}}$ & $\mathtt{\bar{1}}$\\\hline
\end{tabular}
\right)  =3\rightarrow\left(
\begin{tabular}
[c]{|l|}\hline
$\mathtt{3}$\\\hline
$\mathtt{\bar{3}}$\\\hline
\end{tabular}
\left(  1\rightarrow%
\begin{tabular}
[c]{|l|}\hline
$\mathtt{3}$\\\hline
$\mathtt{\bar{1}}$\\\hline
\end{tabular}
\right)  \right)  \\
=3\rightarrow%
\begin{tabular}
[c]{|l|l|l}\hline
$\mathtt{3}$ & $\mathtt{\bar{3}}$ & \multicolumn{1}{|l|}{$\mathtt{\bar{1}}$%
}\\\hline
$\mathtt{\bar{3}}$ & $\mathtt{\bar{1}}$ & \\\cline{1-2}%
\end{tabular}
=%
\begin{tabular}
[c]{|l|}\hline
$\mathtt{2}$\\\hline
$\mathtt{3}$\\\hline
\end{tabular}
\left(  \overline{2}\rightarrow%
\begin{tabular}
[c]{|l|l}\hline
$\mathtt{\bar{3}}$ & \multicolumn{1}{|l|}{$\mathtt{\bar{1}}$}\\\hline
$\mathtt{\bar{1}}$ & \\\cline{1-1}%
\end{tabular}
\right)  =%
\begin{tabular}
[c]{|l|l|ll}\hline
$\mathtt{2}$ & $\mathtt{\bar{3}}$ & $\mathtt{\bar{1}}$ &
\multicolumn{1}{|l|}{$\mathtt{\bar{1}}$}\\\hline
$\mathtt{3}$ & $\mathtt{\bar{2}}$ &  & \\\cline{1-2}%
\end{tabular}
\end{multline*}
\end{example}

Finally for any vertex $w\in\Gamma$, we will have:%
\begin{align*}
P(w)  &  =%
\begin{tabular}
[c]{|l|}\hline
$w$\\\hline
\end{tabular}
\text{ if }w\text{ is a letter,}\\
P(w)  &  =x\rightarrow P(u)\text{ if }w=ux\text{ with }u\text{ a word and
}x\text{ a letter.}%
\end{align*}

\subsection{Robinson Schensted correspondence}

In this section a bijection is established between words $w$ of length $l$ on
$\mathcal{G}\ $and pairs $(P(w),Q(w))$ where $P(w)$ is the tableau of type
$G_{2}$ defined above and $Q(w)$ is an oscillating tableau of type $G_{2}$.

\begin{definition}
\label{def_tab_osci}An oscillating tableau $Q$ of type $G_{2}$ and length $l$
is a sequence $(Q_{1},...,Q_{l})$ of Young diagrams whose columns have height
$1$ or $2$ satisfying for $k=1,...,l$ one of the following assertions:

\begin{enumerate}
\item $Q_{k+1}$ is obtained by adding one box to $Q_{k}.$

\item $Q_{k+1}$ is obtained by deleting one box in $Q_{k}.$

\item $Q_{k+1}=Q_{k}.$

\item $Q_{k+1}$ is obtained from $Q_{k}$ by moving one box from height $2$ to
height $1.$

\item $Q_{k+1}$ is obtained from $Q_{k}$ by moving one box from height $1$ to
height $2$
\end{enumerate}
\end{definition}

Let $w=x_{1}\cdot\cdot\cdot x_{l}\ $be a word. The construction of $P(w)$
involves the construction of the $l$ tableaux of type $G_{2}$ defined by
$P_{i}=P(x_{1}\cdot\cdot\cdot x_{i}),$ $i=1,...,l.$ For $w\in\mathcal{G}%
^{\ast}$ we denote by $Q(w)$ the sequence of shapes of the tableaux
$(P_{1},...,P_{l})$.

\begin{proposition}
\label{prop_Q(w)_oscill}$Q(w)$ is an oscillating tableau of type $G_{2}.$
\end{proposition}

\begin{proof}
It follows immediately from the proof of Lemma \ref{lem_aux}.
\end{proof}

\begin{theorem}
\label{th_Q_symbol}Two vertices $w_{1}$ and $w_{2}$ of $\Gamma$ belong to the
same connected component if and only if $Q(w_{1})=Q(w_{2})$.
\end{theorem}

\begin{proof}
The proof is the same than in Theorem 3.4.3 of \cite{L2(RS)}
\end{proof}

\begin{corollary}
Let $\mathcal{G}_{l}^{\ast}$ and $\mathcal{O}_{l}$ be the set of words of
length $l$ on $\mathcal{G}$ and the set of pairs $(P,Q)$ where $P$ is a
tableau of type $G_{2}$ and $Q$ an oscillating tableau of type $G_{2}$ and
length $l$ such that $P$ has shape $Q_{l}$ ($Q_{l}$ is the last shape of $Q$).
Then the map:%
\[
\mathcal{%
\begin{tabular}
[c]{l}%
$\Psi:\mathcal{G}_{l}^{\ast}\rightarrow\mathcal{O}_{l}$\\
$w\mapsto(P(w),Q(w))$%
\end{tabular}
}%
\]
is a bijection.
\end{corollary}

\begin{proof}
The proof is analogous to that of Theorem 5.2.2 of \cite{L2(RS)}.
\end{proof}

\section{Computing the global basis of $V(\lambda)$}

\subsection{The representation $V(\Lambda_{1})$}

The vector representation $V(\Lambda_{1})$ of $U_{q}(G_{2})$ is the vector
space of basis $\{v_{x},$ $x\in\mathcal{G}\}$ where
\[
t_{i}(v_{x})=q_{i}^{<\mathrm{wt}(x),\alpha_{i}>}\text{ for }i=1,2.
\]%
\begin{equation}
\text{and }\left\{
\begin{tabular}
[c]{l}%
for $x\neq0,$ $f_{i}(v_{x})=v_{y}$ if $\widetilde{f}_{i}(x)=y,$ otherwise
$f_{i}(v_{x})=0$\\
$f_{i}(v_{0})=\delta_{i,1}(q+q^{-1})v_{\overline{3}}$\\
for $x\neq0,$ $e_{i}(v_{x})=v_{y}$ if $\widetilde{e}_{i}(x)=y,$ otherwise
$e_{i}(v_{x})=0$\\
$e_{i}(v_{0})=\delta_{i,1}(q+q^{-1})v_{3}$%
\end{tabular}
\right.  . \label{action fi_vect}%
\end{equation}
Note that, with our definition of the action of $f_{1}$ on $v_{n}$ and
$v_{0},$ we have $f_{1}^{(2)}(v_{n})=v_{\overline{n}}$.

\begin{remark}
\label{remarkL1}Set $L_{1}=\underset{x\in\mathcal{G}}{%
%TCIMACRO{\tbigoplus }%
%BeginExpansion
{\textstyle\bigoplus}
%EndExpansion
}Av_{x}.\;$If we identify the image of $v_{x}$ by the canonical projection
$L_{1}\rightarrow L_{1}/qL_{1}$ with $x$ then $(L_{1},B_{1}=B(\Lambda_{1}))$
is the crystal basis of $V(\Lambda_{1}).$
\end{remark}

\noindent The basis $\{v_{x},$ $x\in\mathcal{G}\}$ satisfies conditions
(\ref{cond_cong}) and (\ref{cond_invo}) of Theorem \ref{Th_K} thus is the
canonical basis of $V(\Lambda_{1}).\;$We have $G(%
\begin{tabular}
[c]{|l|}\hline
$x$\\\hline
\end{tabular}
)=v_{x}$ for any $x\in B(\Lambda_{1}).$

\subsection{Canonical basis of $V(\Lambda_{2})$}

Similarly to \cite{lec3} we introduce $W(\Lambda_{2})$ the $q$-analogue to the
$2$-th wedge product of $V(\Lambda_{1})$ defined by%
\[
W(\Lambda_{2})=V(\Lambda_{1})^{\otimes2}/N
\]
where $N$ is the sub-module isomorphic to $V(2\Lambda_{1})$ appearing in the
decomposition
\[
V(\Lambda_{1})^{\otimes2}\simeq V(2\Lambda_{1})%
%TCIMACRO{\tbigoplus }%
%BeginExpansion
{\textstyle\bigoplus}
%EndExpansion
V(\Lambda_{2})%
%TCIMACRO{\tbigoplus }%
%BeginExpansion
{\textstyle\bigoplus}
%EndExpansion
V(\Lambda_{1})%
%TCIMACRO{\tbigoplus }%
%BeginExpansion
{\textstyle\bigoplus}
%EndExpansion
V(0)
\]
of $V(\Lambda_{1})^{\otimes2}$ into its irreducible components. Let $\Psi$ be
the canonical projection $V(\Lambda_{1})^{\otimes2}\rightarrow W(\Lambda_{2}%
)$.\ Set%
\begin{equation}
\Psi(v_{x_{1}}\otimes v_{x_{2}})=v_{x_{1}}\wedge v_{x_{2}}.\label{def_ext}%
\end{equation}

\begin{proposition}
\label{prop_relation}In $W(\Lambda_{2})$ we have the relations:

\begin{enumerate}
\item  for $x\neq0$, $v_{x}\wedge v_{x}=0,$

\item  for $xy\in B(\Lambda_{2})$ and $x\neq\overline{y}$, $v_{y}\wedge
v_{x}=\left\{
\begin{tabular}
[c]{l}%
$-q^{2}v_{x}\wedge v_{y}$ if $x=0$ or $y=0$\\
$-q^{3}v_{x}\wedge v_{y}$ if $(x,y)=(2,3)$ or $(x,y)=(\overline{3}%
,\overline{2})$\\
$-qv_{x}\wedge v_{y}$ otherwise
\end{tabular}
\right.  ,$

\item
\[
\left\{
\begin{tabular}
[c]{l}%
$v_{0}\wedge v_{1}=-q^{2}v_{1}\wedge v_{0}+(q^{5}-q)v_{2}\wedge v_{3}$\\
$v_{\overline{3}}\wedge v_{1}=-q^{3}v_{1}\wedge v_{\overline{3}}%
+(q^{3}-q)v_{2}\wedge v_{0}$\\
$v_{\overline{2}}\wedge v_{1}=-q^{3}v_{1}\wedge v_{\overline{2}}%
+(q^{3}-q)v_{3}\wedge v_{0}$%
\end{tabular}
\right.  \text{ and }\left\{
\begin{tabular}
[c]{l}%
$v_{\overline{1}}\wedge v_{0}=-q^{2}v_{0}\wedge v_{\overline{1}}%
+(q^{5}-q)v_{\overline{3}}\wedge v_{\overline{2}}$\\
$v_{\overline{1}}\wedge v_{3}=-q^{3}v_{3}\wedge v_{\overline{1}}%
+(q^{3}-q)v_{0}\wedge v_{\overline{2}}$\\
$v_{\overline{1}}\wedge v_{2}=-q^{3}v_{2}\wedge v_{\overline{1}}%
+(q^{3}-q)v_{0}\wedge v_{\overline{3}}$%
\end{tabular}
\right.  ,
\]

\item $v_{\overline{3}}\wedge v_{3}=-q^{4}v_{3}\wedge v_{\overline{3}}%
-qv_{0}\wedge v_{0},$

\item $v_{\overline{2}}\wedge v_{2}=-q^{4}v_{2}\wedge v_{\overline{2}}%
+(q^{7}-q)v_{3}\wedge v_{\overline{3}}+q^{4}v_{0}\wedge v_{0},$

\item $v_{\overline{1}}\wedge v_{1}=-q^{4}v_{1}\wedge v_{\overline{1}}%
+(q^{5}-q^{3})v_{2}\wedge v_{\overline{2}}+(-q^{8}+q^{6}+q^{4}+q^{2}%
)v_{3}\wedge v_{\overline{3}}+(-q^{5}+q^{3}-q)v_{0}\wedge v_{0}.$
\end{enumerate}
\end{proposition}

\begin{proof}
The proof is analogous to that of Proposition 3.1.1 of \cite{lec3}. Starting
from the highest weight vector $v_{1}\otimes v_{1}$ of $N$ we apply
Chevalley's operators $f_{i}$ $i=1,2$ to compute a basis of vectors of
$N.\;$These vectors are annihilated by $\Psi$ which give equalities in
$W(\Lambda_{2}).\;$The relations of the Proposition are then obtained by
linear combinations of these equalities.
\end{proof}

\noindent For any column $C$ of reading $w=c_{1}c_{2}$ where the $c_{i}$'s are
letters, we set $v_{C}=v_{c_{1}}\wedge v_{c_{2}}$.\ Then each vector
$\Psi(v_{x_{1}}\otimes v_{x_{2}})=v_{x_{1}}\wedge v_{x_{2}}$ can be decomposed
into a linear combination of vectors $v_{C}$ by applying from left to right a
relation given in the above proposition.

\begin{lemma}
The vectors of $\{v_{C},$ $C\in\mathbf{C}(2)\}$ form a basis of $W(\Lambda_{2}).$
\end{lemma}

\begin{proof}
Each vector of $W(\Lambda_{2})$ can be decomposed into a linear combination of
vectors $v_{C}.\;$So it suffices to proves that $\dim(W(\Lambda_{2}%
))=\mathrm{card}(\mathbf{C}(2)).$ It follows from the definition of
$W(\Lambda_{2})$ that $\dim(W(\Lambda_{2}))=\dim(V(\Lambda_{1}))^{2}%
-\dim(V(2\Lambda_{1}))=49-27=22$ which is exactly the number of columns of
height $2$ on $\mathcal{G}$.
\end{proof}

The coordinates of a vector $v_{x_{1}}\wedge v_{x_{2}}$ on the basis
$\{v_{C},$ $C\in\mathbf{C}(2)\}$ are all in $\mathbb{Z[}q]$ since it is true
for the coefficients appearing in the relations of Proposition
\ref{prop_relation}.

\noindent Consider the $A$-lattice $L_{2}=\underset{C\in\mathbf{C}(2)}{%
%TCIMACRO{\tbigoplus }%
%BeginExpansion
{\textstyle\bigoplus}
%EndExpansion
}Av_{C}$ of $W(\Lambda_{2})$ and denote by $\pi_{2}$ the projection
$L_{2}\rightarrow L_{2}/qL_{2}$.\ We identify $\pi_{2}(v_{C})$ with the word
$\mathrm{w}(C)$.

\begin{lemma}
\label{lem_base_cryst_Wp}$(L_{2},B_{2}=\{\mathrm{w}(C),$ $C\in\mathbf{C}%
(2)\})$ is a crystal basis of $W(\Lambda_{2})$.
\end{lemma}

\begin{proof}
The proof is the same than in Lemma 3.1.3 of \cite{lec3}.
\end{proof}

The vector $v_{\Lambda_{2}}=v_{1}\wedge v_{2}$ is of highest weight
$\Lambda_{2}$ in $W(\Lambda_{2}).$ We identify $V(\Lambda_{2})$ with the
sub-module of $W(\Lambda_{2})$ isomorphic to $V(\Lambda_{2}).$ In the sequel
we need the explicit description of the action of the Chevalley operators
$f_{1}$ and $f_{2}$ on the basis $\{v_{C},$ $C\in\mathbf{C}(2)\}.$ They are
given by the following tables (where we have written for short $C$ in place of
$v_{C})$ obtained from (\ref{action fi_vect}), (\ref{def_ext}) and Proposition
\ref{prop_relation}:%
\[
{%
\begin{tabular}
[c]{c|c}%
\begin{tabular}
[c]{l}%
$C$\\
\ \
\end{tabular}
&
\begin{tabular}
[c]{l}%
$f_{1}(v_{C})$\\
\ \
\end{tabular}
\\\hline%
\begin{tabular}
[c]{l}%
\ \ \\\hline
\multicolumn{1}{|l|}{$\mathtt{1}$}\\\hline
\multicolumn{1}{|l|}{$\mathtt{3}$}\\\hline
\ \
\end{tabular}
& $%
\begin{tabular}
[c]{|l|}\hline
$\mathtt{2}$\\\hline
$\mathtt{3}$\\\hline
\end{tabular}
+q%
\begin{tabular}
[c]{|l|}\hline
$\mathtt{1}$\\\hline
$\mathtt{0}$\\\hline
\end{tabular}
$\\\hline%
\begin{tabular}
[c]{l}%
\ \ \\\hline
\multicolumn{1}{|l|}{$\mathtt{1}$}\\\hline
\multicolumn{1}{|l|}{$\mathtt{0}$}\\\hline
\ \
\end{tabular}
& $%
\begin{tabular}
[c]{|l|}\hline
$\mathtt{2}$\\\hline
$\mathtt{0}$\\\hline
\end{tabular}
+(q^{2}+1)%
\begin{tabular}
[c]{|l|}\hline
$\mathtt{1}$\\\hline
$\mathtt{\bar{3}}$\\\hline
\end{tabular}
$\\\hline%
\begin{tabular}
[c]{l}%
\ \ \\\hline
\multicolumn{1}{|l|}{$\mathtt{1}$}\\\hline
\multicolumn{1}{|l|}{$\mathtt{\bar{3}}$}\\\hline
\ \
\end{tabular}
&
\begin{tabular}
[c]{|l|}\hline
$\mathtt{2}$\\\hline
$\mathtt{\bar{3}}$\\\hline
\end{tabular}
\\\hline%
\begin{tabular}
[c]{l}%
\ \ \\\hline
\multicolumn{1}{|l|}{$\mathtt{1}$}\\\hline
\multicolumn{1}{|l|}{$\mathtt{\bar{2}}$}\\\hline
\ \
\end{tabular}
& $%
\begin{tabular}
[c]{|l|}\hline
$\mathtt{2}$\\\hline
$\mathtt{\bar{2}}$\\\hline
\end{tabular}
+q%
\begin{tabular}
[c]{|l|}\hline
$\mathtt{1}$\\\hline
$\mathtt{\bar{1}}$\\\hline
\end{tabular}
$\\\hline%
\begin{tabular}
[c]{l}%
\ \ \\\hline
\multicolumn{1}{|l|}{$\mathtt{1}$}\\\hline
\multicolumn{1}{|l|}{$\mathtt{\bar{1}}$}\\\hline
\ \
\end{tabular}
&
\begin{tabular}
[c]{|l|}\hline
$\mathtt{2}$\\\hline
$\mathtt{\bar{1}}$\\\hline
\end{tabular}
\\\hline%
\begin{tabular}
[c]{l}%
\ \ \\\hline
\multicolumn{1}{|l|}{$\mathtt{2}$}\\\hline
\multicolumn{1}{|l|}{$\mathtt{3}$}\\\hline
\ \
\end{tabular}
& $q^{-1}%
\begin{tabular}
[c]{|l|}\hline
$\mathtt{2}$\\\hline
$\mathtt{0}$\\\hline
\end{tabular}
$\\\hline%
\begin{tabular}
[c]{l}%
\ \ \\\hline
\multicolumn{1}{|l|}{$\mathtt{2}$}\\\hline
\multicolumn{1}{|l|}{$\mathtt{0}$}\\\hline
\ \
\end{tabular}
& $(1+q^{-2})%
\begin{tabular}
[c]{|l|}\hline
$\mathtt{2}$\\\hline
$\mathtt{\bar{3}}$\\\hline
\end{tabular}
$\\\hline%
\begin{tabular}
[c]{l}%
\ \ \\\hline
\multicolumn{1}{|l|}{$\mathtt{2}$}\\\hline
\multicolumn{1}{|l|}{$\mathtt{\bar{2}}$}\\\hline
\ \
\end{tabular}
& $q^{-1}%
\begin{tabular}
[c]{|l|}\hline
$\mathtt{2}$\\\hline
$\mathtt{\bar{1}}$\\\hline
\end{tabular}
$%
\end{tabular}
}\ \ \ \ \ \ \ \ \ \ {%
\begin{tabular}
[c]{c|c}%
\begin{tabular}
[c]{l}%
$C$\\
\ \
\end{tabular}
&
\begin{tabular}
[c]{l}%
$f_{1}(v_{C})$\\
\ \
\end{tabular}
\\\hline%
\begin{tabular}
[c]{l}%
\ \ \\\hline
\multicolumn{1}{|l|}{$\mathtt{3}$}\\\hline
\multicolumn{1}{|l|}{$\mathtt{0}$}\\\hline
\ \
\end{tabular}
& $%
\begin{tabular}
[c]{|l|}\hline
$\mathtt{0}$\\\hline
$\mathtt{0}$\\\hline
\end{tabular}
+(q^{3}+q)$%
\begin{tabular}
[c]{|l|}\hline
$\mathtt{3}$\\\hline
$\mathtt{\bar{3}}$\\\hline
\end{tabular}
\\\hline%
\begin{tabular}
[c]{l}%
\ \ \\\hline
\multicolumn{1}{|l|}{$\mathtt{3}$}\\\hline
\multicolumn{1}{|l|}{$\mathtt{\bar{3}}$}\\\hline
\ \
\end{tabular}
&
\begin{tabular}
[c]{|l|}\hline
$\mathtt{0}$\\\hline
$\mathtt{\bar{3}}$\\\hline
\end{tabular}
\\\hline%
\begin{tabular}
[c]{l}%
\ \ \\\hline
\multicolumn{1}{|l|}{$\mathtt{3}$}\\\hline
\multicolumn{1}{|l|}{$\mathtt{\bar{2}}$}\\\hline
\ \
\end{tabular}
& $%
\begin{tabular}
[c]{|l|}\hline
$\mathtt{0}$\\\hline
$\mathtt{\bar{2}}$\\\hline
\end{tabular}
+q^{2}%
\begin{tabular}
[c]{|l|}\hline
$\mathtt{3}$\\\hline
$\mathtt{\bar{1}}$\\\hline
\end{tabular}
$\\\hline%
\begin{tabular}
[c]{l}%
\ \ \\\hline
\multicolumn{1}{|l|}{$\mathtt{3}$}\\\hline
\multicolumn{1}{|l|}{$\mathtt{\bar{1}}$}\\\hline
\ \
\end{tabular}
& $%
\begin{tabular}
[c]{|l|}\hline
$\mathtt{0}$\\\hline
$\mathtt{\bar{1}}$\\\hline
\end{tabular}
$\\\hline%
\begin{tabular}
[c]{l}%
\ \ \\\hline
\multicolumn{1}{|l|}{$\mathtt{0}$}\\\hline
\multicolumn{1}{|l|}{$\mathtt{0}$}\\\hline
\ \
\end{tabular}
& $(1-q^{2})(q+q^{-1})%
\begin{tabular}
[c]{|l|}\hline
$\mathtt{0}$\\\hline
$\mathtt{\bar{3}}$\\\hline
\end{tabular}
$\\\hline%
\begin{tabular}
[c]{l}%
\ \ \\\hline
\multicolumn{1}{|l|}{$\mathtt{0}$}\\\hline
\multicolumn{1}{|l|}{$\mathtt{\bar{2}}$}\\\hline
\ \
\end{tabular}
& $(q+q^{-1})%
\begin{tabular}
[c]{|l|}\hline
$\mathtt{\bar{3}}$\\\hline
$\mathtt{\bar{2}}$\\\hline
\end{tabular}
+%
\begin{tabular}
[c]{|l|}\hline
$\mathtt{0}$\\\hline
$\mathtt{\bar{1}}$\\\hline
\end{tabular}
$\\\hline%
\begin{tabular}
[c]{l}%
\ \ \\\hline
\multicolumn{1}{|l|}{$\mathtt{0}$}\\\hline
\multicolumn{1}{|l|}{$\mathtt{\bar{1}}$}\\\hline
\ \
\end{tabular}
& $(q+q^{-1})%
\begin{tabular}
[c]{|l|}\hline
$\mathtt{\bar{3}}$\\\hline
$\mathtt{\bar{1}}$\\\hline
\end{tabular}
$\\\hline%
\begin{tabular}
[c]{l}%
\ \ \\\hline
\multicolumn{1}{|l|}{$\mathtt{\bar{3}}$}\\\hline
\multicolumn{1}{|l|}{$\mathtt{\bar{2}}$}\\\hline
\ \
\end{tabular}
& $q^{-2}%
\begin{tabular}
[c]{|l|}\hline
$\mathtt{\bar{3}}$\\\hline
$\mathtt{\bar{1}}$\\\hline
\end{tabular}
$%
\end{tabular}
}\text{ \ \ \ and }f_{1}(v_{C})=0\text{ otherwise.}%
\]%

\[
{%
\begin{tabular}
[c]{c|c}%
$%
\begin{tabular}
[c]{l}%
$C$\\
\ \
\end{tabular}
$ & $%
\begin{tabular}
[c]{l}%
$f_{2}(v_{C})$\\
\ \
\end{tabular}
$\\\hline%
\begin{tabular}
[c]{l}%
\ \ \\\hline
\multicolumn{1}{|l|}{$\mathtt{1}$}\\\hline
\multicolumn{1}{|l|}{$\mathtt{2}$}\\\hline
\ \
\end{tabular}
&
\begin{tabular}
[c]{|l|}\hline
$\mathtt{1}$\\\hline
$\mathtt{3}$\\\hline
\end{tabular}
\\\hline%
\begin{tabular}
[c]{l}%
\ \ \\\hline
\multicolumn{1}{|l|}{$\mathtt{2}$}\\\hline
\multicolumn{1}{|l|}{$x$}\\\hline
\ \
\end{tabular}
with $x\in\{\bar{1},0\}$ & $%
\begin{tabular}
[c]{|l|}\hline
$\mathtt{3}$\\\hline
$x$\\\hline
\end{tabular}
$\\\hline%
\begin{tabular}
[c]{l}%
\ \ \\\hline
\multicolumn{1}{|l|}{$\mathtt{\bar{3}}$}\\\hline
\multicolumn{1}{|l|}{$\mathtt{\bar{1}}$}\\\hline
\ \
\end{tabular}
&
\begin{tabular}
[c]{|l|}\hline
$\mathtt{\bar{2}}$\\\hline
$\mathtt{\bar{1}}$\\\hline
\end{tabular}
\\\hline%
\begin{tabular}
[c]{l}%
\ \ \\\hline
\multicolumn{1}{|l|}{$x$}\\\hline
\multicolumn{1}{|l|}{$\mathtt{\bar{3}}$}\\\hline
\ \
\end{tabular}
with $x\in\{0,1\}$ & $%
\begin{tabular}
[c]{|l|}\hline
$x$\\\hline
$\mathtt{\bar{2}}$\\\hline
\end{tabular}
$\\\hline%
\begin{tabular}
[c]{l}%
\ \ \\\hline
\multicolumn{1}{|l|}{$\mathtt{2}$}\\\hline
\multicolumn{1}{|l|}{$\mathtt{\bar{3}}$}\\\hline
\ \
\end{tabular}
& $%
\begin{tabular}
[c]{|l|}\hline
$\mathtt{3}$\\\hline
$\mathtt{\bar{3}}$\\\hline
\end{tabular}
+q^{3}%
\begin{tabular}
[c]{|l|}\hline
$\mathtt{2}$\\\hline
$\mathtt{\bar{2}}$\\\hline
\end{tabular}
$\\\hline%
\begin{tabular}
[c]{l}%
\ \ \\\hline
\multicolumn{1}{|l|}{$\mathtt{2}$}\\\hline
\multicolumn{1}{|l|}{$\mathtt{\bar{2}}$}\\\hline
\ \
\end{tabular}
& $%
\begin{tabular}
[c]{|l|}\hline
$\mathtt{3}$\\\hline
$\mathtt{\bar{2}}$\\\hline
\end{tabular}
$\\\hline%
\begin{tabular}
[c]{l}%
\ \ \\\hline
\multicolumn{1}{|l|}{$\mathtt{3}$}\\\hline
\multicolumn{1}{|l|}{$\mathtt{\bar{3}}$}\\\hline
\ \
\end{tabular}
& $q^{-3}%
\begin{tabular}
[c]{|l|}\hline
$\mathtt{3}$\\\hline
$\mathtt{\bar{2}}$\\\hline
\end{tabular}
$%
\end{tabular}
}\text{ \ \ and }f_{2}(v_{C})=0\text{ otherwise.}%
\]

Now we are going to give the explicit decomposition of the canonical basis of
$V(\Lambda_{2})$ on $\{v_{C},$ $C\in\mathbf{C}(2)\}.\;$Consider an admissible
column $C$ of height $2.\;$If $\mathrm{w}(C)\neq0\overline{2}$ there is a
unique path in $B(\Lambda_{2})$ joining $12$ to $\mathrm{w}(C)$.\ Otherwise we
choose the path $0\overline{2}=\widetilde{f}_{1}\widetilde{f}_{2}%
^{2}\widetilde{f}_{1}^{3}\widetilde{f}_{2}(12).\;$Then (with our choice of the
path joining $12$ to $0\overline{2})$ we can write $\mathrm{w}(C)=\widetilde
{f}_{i_{1}}^{p_{1}}\cdot\cdot\cdot\widetilde{f}_{i_{r}}^{p_{r}}(12)$ with
$i_{k}\neq i_{k+1}$ for $k=1,...,r-1.$

\begin{theorem}
For any admissible column $C$ of height $2,$ $G(C)=f_{i_{1}}^{(p_{1})}%
\cdot\cdot\cdot f_{i_{r}}^{(p_{r})}(v_{\Lambda_{2}})$.
\end{theorem}

\begin{proof}
The vectors $f_{i_{1}}^{(p_{1})}\cdot\cdot\cdot f_{i_{r}}^{(p_{r})}%
(v_{\Lambda_{2}})$ belong to $V_{\mathbb{Q}}(\Lambda_{2})$ and are fixed by
the involution $\overline{\text{%
\begin{tabular}
[c]{l}%
\ \
\end{tabular}
}}.$ So it suffices to prove that the coordinates of the decomposition of each
vector $f_{i_{1}}^{(p_{1})}\cdot\cdot\cdot f_{i_{r}}^{(p_{r})}(v_{\Lambda_{2}%
})$ on the basis $\{v_{C},$ $C\in\mathbf{C}(2)\}$ are all in $\mathbb{Z[}q]$
and such that
\[
f_{i_{1}}^{(p_{1})}\cdot\cdot\cdot f_{i_{r}}^{(p_{r})}(v_{\Lambda_{2}%
})=\mathrm{w}(C)\text{ }\operatorname{mod}(q).
\]
This is shown by an explicit computation from the action of the operators
$f_{1}$ and $f_{2}$ given above.\ The results are given in the table below.
\end{proof}%

\[
{%
\begin{tabular}
[c]{c|c}%
$%
\begin{tabular}
[c]{l}%
$C$\\
\ \
\end{tabular}
$ & $%
\begin{tabular}
[c]{l}%
$G(C)$\\
\ \
\end{tabular}
$\\\hline%
\begin{tabular}
[c]{l}%
\ \ \\\hline
\multicolumn{1}{|l|}{$\mathtt{1}$}\\\hline
\multicolumn{1}{|l|}{$\mathtt{2}$}\\\hline
\ \
\end{tabular}
& $%
\begin{tabular}
[c]{|l|}\hline
$\mathtt{1}$\\\hline
$\mathtt{2}$\\\hline
\end{tabular}
$\\\hline%
\begin{tabular}
[c]{l}%
\ \ \\\hline
\multicolumn{1}{|l|}{$\mathtt{1}$}\\\hline
\multicolumn{1}{|l|}{$\mathtt{3}$}\\\hline
\ \
\end{tabular}
& $%
\begin{tabular}
[c]{|l|}\hline
$\mathtt{1}$\\\hline
$\mathtt{3}$\\\hline
\end{tabular}
$\\\hline%
\begin{tabular}
[c]{l}%
\ \ \\\hline
\multicolumn{1}{|l|}{$\mathtt{2}$}\\\hline
\multicolumn{1}{|l|}{$\mathtt{3}$}\\\hline
\ \
\end{tabular}
& $%
\begin{tabular}
[c]{|l|}\hline
$\mathtt{2}$\\\hline
$\mathtt{3}$\\\hline
\end{tabular}
+q%
\begin{tabular}
[c]{|l|}\hline
$\mathtt{1}$\\\hline
$\mathtt{0}$\\\hline
\end{tabular}
$\\\hline%
\begin{tabular}
[c]{l}%
\ \ \\\hline
\multicolumn{1}{|l|}{$\mathtt{2}$}\\\hline
\multicolumn{1}{|l|}{$\mathtt{\bar{3}}$}\\\hline
\ \
\end{tabular}
& $%
\begin{tabular}
[c]{|l|}\hline
$\mathtt{2}$\\\hline
$\mathtt{\bar{3}}$\\\hline
\end{tabular}
$\\\hline%
\begin{tabular}
[c]{l}%
\ \ \\\hline
\multicolumn{1}{|l|}{$\mathtt{2}$}\\\hline
\multicolumn{1}{|l|}{$\mathtt{0}$}\\\hline
\ \
\end{tabular}
& $%
\begin{tabular}
[c]{|l|}\hline
$\mathtt{2}$\\\hline
$\mathtt{0}$\\\hline
\end{tabular}
+q^{2}%
\begin{tabular}
[c]{|l|}\hline
$\mathtt{1}$\\\hline
$\mathtt{\bar{3}}$\\\hline
\end{tabular}
$\\\hline%
\begin{tabular}
[c]{l}%
\ \ \\\hline
\multicolumn{1}{|l|}{$\mathtt{3}$}\\\hline
\multicolumn{1}{|l|}{$\mathtt{0}$}\\\hline
\ \
\end{tabular}
& $%
\begin{tabular}
[c]{|l|}\hline
$\mathtt{3}$\\\hline
$\mathtt{0}$\\\hline
\end{tabular}
+q^{2}%
\begin{tabular}
[c]{|l|}\hline
$\mathtt{1}$\\\hline
$\mathtt{\bar{2}}$\\\hline
\end{tabular}
$\\\hline%
\begin{tabular}
[c]{l}%
\ \ \\\hline
\multicolumn{1}{|l|}{$\mathtt{3}$}\\\hline
\multicolumn{1}{|l|}{$\mathtt{\bar{3}}$}\\\hline
\ \
\end{tabular}
& $%
\begin{tabular}
[c]{|l|}\hline
$\mathtt{3}$\\\hline
$\mathtt{\bar{3}}$\\\hline
\end{tabular}
+q^{3}%
\begin{tabular}
[c]{|l|}\hline
$\mathtt{3}$\\\hline
$\mathtt{\bar{2}}$\\\hline
\end{tabular}
$%
\end{tabular}
\ \ \ \ \ \ \ \
\begin{tabular}
[c]{c|c}%
$%
\begin{tabular}
[c]{l}%
$C$\\
\ \
\end{tabular}
$ & $%
\begin{tabular}
[c]{l}%
$G(C)$\\
\ \
\end{tabular}
$\\\hline%
\begin{tabular}
[c]{l}%
\ \ \\\hline
\multicolumn{1}{|l|}{$\mathtt{\bar{2}}$}\\\hline
\multicolumn{1}{|l|}{$\mathtt{\bar{1}}$}\\\hline
\ \
\end{tabular}
& $%
\begin{tabular}
[c]{|l|}\hline
$\mathtt{\bar{2}}$\\\hline
$\mathtt{\bar{1}}$\\\hline
\end{tabular}
$\\\hline%
\begin{tabular}
[c]{l}%
\ \ \\\hline
\multicolumn{1}{|l|}{$\mathtt{\bar{3}}$}\\\hline
\multicolumn{1}{|l|}{$\mathtt{\bar{1}}$}\\\hline
\ \
\end{tabular}
& $%
\begin{tabular}
[c]{|l|}\hline
$\mathtt{\bar{3}}$\\\hline
$\mathtt{\bar{1}}$\\\hline
\end{tabular}
$\\\hline%
\begin{tabular}
[c]{l}%
\ \ \\\hline
\multicolumn{1}{|l|}{$\mathtt{\bar{3}}$}\\\hline
\multicolumn{1}{|l|}{$\mathtt{\bar{2}}$}\\\hline
\ \
\end{tabular}
& $%
\begin{tabular}
[c]{|l|}\hline
$\mathtt{\bar{3}}$\\\hline
$\mathtt{\bar{2}}$\\\hline
\end{tabular}
+q%
\begin{tabular}
[c]{|l|}\hline
$\mathtt{0}$\\\hline
$\mathtt{\bar{1}}$\\\hline
\end{tabular}
$\\\hline%
\begin{tabular}
[c]{l}%
\ \ \\\hline
\multicolumn{1}{|l|}{$\mathtt{3}$}\\\hline
\multicolumn{1}{|l|}{$\mathtt{\bar{2}}$}\\\hline
\ \
\end{tabular}
& $%
\begin{tabular}
[c]{|l|}\hline
$\mathtt{3}$\\\hline
$\mathtt{\bar{2}}$\\\hline
\end{tabular}
$\\\hline%
\begin{tabular}
[c]{l}%
\ \ \\\hline
\multicolumn{1}{|l|}{$\mathtt{0}$}\\\hline
\multicolumn{1}{|l|}{$\mathtt{\bar{2}}$}\\\hline
\ \
\end{tabular}
& $%
\begin{tabular}
[c]{|l|}\hline
$\mathtt{0}$\\\hline
$\mathtt{\bar{2}}$\\\hline
\end{tabular}
+q^{2}%
\begin{tabular}
[c]{|l|}\hline
$\mathtt{3}$\\\hline
$\mathtt{\bar{1}}$\\\hline
\end{tabular}
$\\\hline%
\begin{tabular}
[c]{l}%
\ \ \\\hline
\multicolumn{1}{|l|}{$\mathtt{0}$}\\\hline
\multicolumn{1}{|l|}{$\mathtt{\bar{3}}$}\\\hline
\ \
\end{tabular}
& $%
\begin{tabular}
[c]{|l|}\hline
$\mathtt{0}$\\\hline
$\mathtt{\bar{3}}$\\\hline
\end{tabular}
+q^{2}%
\begin{tabular}
[c]{|l|}\hline
$\mathtt{2}$\\\hline
$\mathtt{\bar{1}}$\\\hline
\end{tabular}
$\\\hline%
\begin{tabular}
[c]{l}%
\ \ \\\hline
\multicolumn{1}{|l|}{$\mathtt{0}$}\\\hline
\multicolumn{1}{|l|}{$\mathtt{0}$}\\\hline
\ \
\end{tabular}
& $%
\begin{tabular}
[c]{|l|}\hline
$\mathtt{0}$\\\hline
$\mathtt{0}$\\\hline
\end{tabular}
+(q^{3}+q)%
\begin{tabular}
[c]{|l|}\hline
$\mathtt{3}$\\\hline
$\mathtt{\bar{3}}$\\\hline
\end{tabular}
+q^{2}%
\begin{tabular}
[c]{|l|}\hline
$\mathtt{2}$\\\hline
$\mathtt{\bar{2}}$\\\hline
\end{tabular}
+q^{3}%
\begin{tabular}
[c]{|l|}\hline
$\mathtt{1}$\\\hline
$\mathtt{\bar{1}}$\\\hline
\end{tabular}
$%
\end{tabular}
}%
\]

\subsection{Algorithm for the global basis of $V(\lambda).$}

\subsubsection{The representation $W(\lambda)$}

To make our notation homogeneous write $W(\Lambda_{1})=V(\Lambda_{1})$.

\noindent Consider $\lambda=\lambda_{1}\Lambda_{1}+\lambda_{2}\Lambda_{2}\in
P_{+}$ and set%
\begin{equation}
W(\lambda)=W(\Lambda_{1})^{\otimes\lambda_{1}}\otimes W(\Lambda_{2}%
)^{\otimes\lambda_{2}}. \label{W(lambda)}%
\end{equation}
The natural basis of $W(\lambda)$ consists of the tensor products $v_{C_{r}%
}\otimes\cdot\cdot\cdot\otimes v_{C_{1}}$ of basis vectors $v_{C}$ of the
previous section appearing in (\ref{W(lambda)}).\ The juxtaposition of the
columns $C_{1},...,C_{r}$ is called a tabloid of shape $Y(\lambda)$.\ We can
regard it as a filling $\tau$ of the Young diagram of shape $\lambda$ the
$i$-th column of which is equal to $C_{i}$. We shall write $v_{\tau}=v_{C_{r}%
}\otimes\cdot\cdot\cdot\otimes v_{C_{1}}.$ Note that the columns of $\tau$ are
not necessarily admissible and there is no condition on the rows. The reading
of the tabloid $\tau=C_{1}\cdot\cdot\cdot C_{r}$ is $\mathrm{w(}%
\tau)=\mathrm{w(}C_{r})\cdot\cdot\cdot\mathrm{w(}C_{1}).$ We denote by
$\mathbf{T}(\lambda)$ the set of tabloids of shape $\lambda.$

Let $L_{\lambda}$ be the $A$-submodule of $W(\lambda)$ generated by the
vectors $v_{\tau},$ $\tau\in\mathbf{T}(\lambda)$. We identify the image of the
vector $v_{\tau}$ by the projection $\pi_{\lambda}:L_{\lambda}\rightarrow
L_{\lambda}/qL_{\lambda}$ with the word \textrm{w(}$\tau)$. The pair
$(L_{\lambda},B_{\lambda}=\{\mathrm{w(}\tau),$ $\tau\in\mathbf{T}(\lambda)\})$
is then a crystal basis of $W(\lambda)$. Indeed by Remark \ref{remarkL1} and
Lemma \ref{lem_base_cryst_Wp}, it is the tensor product of the crystal bases
of the representations $W(\Lambda_{p})$ $p=1,2$ occurring in $W(\lambda)$. Set%
\[
v_{\lambda}=v_{\Lambda_{1}}^{\otimes\lambda_{1}}\otimes v_{\Lambda_{2}%
}^{\otimes}.
\]
We identify $V(\lambda)$ with the submodule of $W(\lambda)$ generated by
$v_{\lambda}$. Then, with the above notations, $v_{\lambda}=v_{T_{\lambda}}$
where $T_{\lambda}$ is the orthogonal tableau of shape $\lambda$ whose $k$-th
row is filled by letters $k$ for $k=1,2$.\ By Theorem 4.2 of \cite{Ka2}, we
know that
\[
B(\lambda)=\{\widetilde{f}_{i_{1}}^{a_{1}}\cdot\cdot\cdot\widetilde{f}_{i_{r}%
}^{a_{r}}\mathrm{w}(T_{\lambda});\text{ }i_{1},...,i_{r}=1,2;\text{ }%
a_{1},...,a_{r}>0\}-\{0\}.
\]
The actions of $\widetilde{e}_{i}$ and $\widetilde{f}_{i},$ $i=1,2$ on each
$\mathrm{w(}\tau)\in B_{\lambda}$ are identical to those obtained by
considering $\mathrm{w(}\tau)$ as a vertex of $\Gamma$ since it is true on
$B_{1}$ and $B_{2}$. Hence, by Theorem \ref{KM}, $B(\lambda)=\{\mathrm{w}(T);$
$T\in\mathbf{T}_{G}(\lambda)\}$. For each vector $G(T)$ of the canonical basis
of $V(\lambda)$ we will have:%
\[
G(T)\equiv v_{T}\operatorname{mod}qL_{\lambda}.
\]
The aim of this section is to describe an algorithm computing the
decomposition of the canonical basis $\{G(T);T\in\mathbf{T}_{G}(\lambda)\}$
onto the basis $\{v_{\tau};\tau\in\mathbf{T}(\lambda)\}$ of $W(\lambda)$.

In the sequel we will need a total order on the readings of the tabloids. Let
$w_{1}=x_{1}\cdot\cdot\cdot x_{l}$ and $w_{2}=y_{1}\cdot\cdot\cdot y_{l}$ be
two distinct vertices of $\Gamma$ with the same length and $k$ the lowest
integer such that $x_{k}\neq y_{k}.$ We write $w_{1}\trianglelefteq w_{2}$ if
$x_{k}\preceq y_{k}$ in $\mathcal{G}$ and $w_{1}\vartriangleright w_{2}$
otherwise that is, $\trianglelefteq$ is the lexicographic order on
$\mathcal{G}^{\ast}$. We endow the set $\mathbf{T}(\lambda)$ with the total
order:%
\[
\tau_{1}\trianglelefteq\tau_{2}\Longleftrightarrow\mathrm{w(}\tau
_{1})\trianglelefteq\mathrm{w(}\tau_{2}).
\]
In fact $\trianglelefteq$ is a lexicographic order defined on the readings of
the tabloids of $\mathbf{T}(\lambda)$ such that $x\trianglelefteq\widetilde
{f}_{i}(x)$ for any letter $x$ with $\varphi_{i}(x)\neq0$.

We are going to compute the canonical basis $\{G(T);T\in\mathbf{T}_{G}%
(\lambda)\}$ in two steps. First we obtain an intermediate basis
$\{A(T);T\in\mathbf{T}_{G}(\lambda)\}$ which is fixed by the involution
$\overline{\text{%
\begin{tabular}
[c]{l}%
\ \
\end{tabular}
}}$ (condition (\ref{cond_invo})). When $T=C$ is an admissible column,
$A(C)=G(C)$.\ In the general case, we have to correct $A(T)$ in order to
verify condition (\ref{cond_cong}). This second step is easy because we can
prove that the transition matrix from $\{A(T);T\in\mathbf{T}_{G}(\lambda)\}$
to $\{G(T);T\in\mathbf{T}_{G}(\lambda)\}$ is unitriangular once the orthogonal
tableaux and the tabloids are ordered by $\trianglelefteq$.

We will need the following lemma whose proof is the same than in Lemma 4.1.1
of \cite{lec3}.

\begin{lemma}
\label{Lem_F()}Let $v\in V(\lambda)$ be a vector of the type
\begin{equation}
v=f_{i_{1}}^{(r_{1})}\cdot\cdot\cdot f_{i_{s}}^{(r_{s})}v_{\lambda
}\label{Type_F()}%
\end{equation}
where $(i_{1},...,i_{s})$ and $(r_{1},...,r_{s})$ are two sequences of
integers.\ Then the coordinates of $v$ on the basis $\{v_{\tau};\tau
\in\mathbf{T}(\lambda)\}$ belong to $\mathbb{Z}[q,q^{-1}]$.
\end{lemma}

\subsubsection{The basis $A(T)$}

The basis $\{A(T)\}$ will be a monomial basis, that is, a basis of the form
\begin{equation}
A(T)=f_{i_{1}}^{(r_{1})}\cdot\cdot\cdot f_{i_{m}}^{(r_{m})}v_{\lambda}.
\label{A(T)_monom}%
\end{equation}
By Lemma \ref{Lem_F()}, the coordinates of $A(T)$ on the basis $\{v_{\tau}\}$
of $W(\lambda)$ belongs to $\mathbb{Z[}q,q^{-1}]$. To find the two sequences
of integers $(i_{1},...,i_{m})$ and $(r_{1},...,r_{m})$ associated to $T$, we
proceed as follows.

\noindent Write $T=C_{1}\cdot\cdot\cdot C_{s}\neq T_{\lambda}\in\mathbf{T}%
_{G}(\lambda)$.\ Let $C_{k}$ be the rightmost column of $T$ such that
$\mathrm{w(}C_{k})$ is not a highest weight vertex (i.e. $\mathrm{w(}%
C_{k})\neq1,$ $12)$.\ When $\mathrm{w(}C_{k})\neq0\overline{2}$ let $i_{1}%
\in\{1,2\}$ be the unique integer such that $\widetilde{e}_{i_{1}}%
(\mathrm{w}(C_{k}))\neq0$. When $\mathrm{w(}C_{k})\neq0\overline{2}$ we choose
$i_{1}=1.$\ If $k=1,$ set $l=1$. If $k>1$ and $f_{i_{1}}v_{C_{k}}\not =0$ or
$\widetilde{e}_{i_{1}}(\mathrm{w}(C_{k-1}))=0$ set $l=k.\;$Otherwise let $l$
be the lowest integer $l<k$ satisfying the two conditions%
\begin{equation}
\left\{
\begin{tabular}
[c]{l}%
$\mathrm{(i):}$ $f_{i_{1}}v_{C_{j}}=0$ for $j=l+1,...,k$\\
$\mathrm{(ii):}$ $\widetilde{e}_{i_{1}}(\mathrm{w}(C_{j}))\neq0$ for
$j=l,...,k$%
\end{tabular}
\right.  . \label{cond}%
\end{equation}
Set $\varepsilon_{i_{1},j}=\varepsilon_{i_{1}}((\mathrm{w}(C_{j}))$ for
$j=l,...,k$ and $r_{1}=\overset{k}{\underset{j=l}{\sum}}\varepsilon_{i_{1}j}%
$.\ Write $T_{1}$ for the tabloid obtained by changing in $T$ each column
$C_{j},$ $j=l,...,k$ into the column of reading $\widetilde{e}_{i_{1}%
}^{\varepsilon_{i_{1},j}}(\mathrm{w}(C_{j})).$

\begin{lemma}
\label{T1_Tab}$T_{1}\in\mathbf{T}_{G}(\lambda).$
\end{lemma}

\begin{proof}
Set $T_{1}=D_{1}\cdot\cdot\cdot D_{s}$.\ Then by construction of $T_{1},$ we
have $D_{i}=C_{i}$ for $i\notin\{l,...,k\}.\;$The columns $C_{i}$ with
$i=k+1,...,s$ are of highest weight.\ One has $\mathrm{w}(D_{k}\cdot\cdot\cdot
D_{s})=\widetilde{e}_{i_{1}}^{\varepsilon_{i_{1},k}}(\mathrm{w(}C_{k}%
\cdot\cdot\cdot C_{s}))$ by (\ref{TENS2}) because $\varepsilon_{i_{1}%
,k}=\varepsilon_{i_{1}}\mathrm{w}(C_{k}\cdot\cdot\cdot C_{s})$ and
$\varphi_{i_{1}}\mathrm{w(}C_{j})=0$ for $j=k+1,...,s$ (otherwise we would
have $\widetilde{e}_{i_{1}}(\mathrm{w(}C_{k}))=0$ since the letters
$1,...,i_{1}$ would occur in $C_{k}$).\ So $D_{k}\cdot\cdot\cdot D_{s}$ is a
tableau of type $G_{2}$.

We have $\mathrm{w(}D_{l}\cdot\cdot\cdot D_{k})=\widetilde{e}_{i_{1}}^{r_{1}%
}\mathrm{w(}(C_{l}\cdot\cdot\cdot C_{k}))$ because $\varphi_{i_{1}}%
\mathrm{w(}C_{j})=0$ for $j=l+1,...,k$ (from \textrm{(i) }of (\ref{cond})) so
is the reading of a tableau of type $G_{2}$. Hence $D_{l}\cdot\cdot\cdot
D_{k}$ is a tableau of type $G_{2}$.

The proof will be complete if we show that $D_{1}\cdot\cdot\cdot D_{l}%
=C_{1}\cdot\cdot\cdot C_{l-1}D_{l}$ is a tableau of type $G_{2},$ that is if
we prove that $C_{l-1}\preceq D_{l}$. We can suppose $l>1.$ Then we have
either $f_{i_{1}}v_{C_{l}}\neq0$ and $\widetilde{e}_{i_{1}}(\mathrm{w}%
(C_{l}))\neq0$ or $\widetilde{e}_{i_{1}}(\mathrm{w}(C_{l-1}))=0.$

\noindent Suppose $\widetilde{e}_{i_{1}}(\mathrm{w}(C_{l-1}))=0.\;$Then we
have $\widetilde{e}_{i_{1}}^{\varepsilon_{i_{1},l}}(\mathrm{w}(C_{l}%
)\mathrm{w}(C_{l-1}))=\mathrm{w}(D_{l})\mathrm{w}(C_{l-1})$ by (\ref{TENS2}%
).\ So $C_{l-1}\preceq D_{l}$.

\noindent Now suppose $f_{i_{1}}v_{C_{l}}\neq0$ and $\widetilde{e}_{i_{1}%
}(\mathrm{w}(C_{l}))\neq0.$ If $i_{1}=2,$ $C_{l}$ is necessarily the column of
reading $3\overline{3}$. We have%
\[
C_{l-1}C_{l}=%
\begin{tabular}
[c]{|l|l|}\hline
$x$ & $\mathtt{3}$\\\hline
$y$ & $\mathtt{\bar{3}}$\\\hline
\end{tabular}
\text{ and }C_{l-1}D_{l}=%
\begin{tabular}
[c]{|l|l|}\hline
$x$ & $\mathtt{2}$\\\hline
$y$ & $\mathtt{\bar{3}}$\\\hline
\end{tabular}
\]
with $x\preceq3$ and $\mathrm{dist}(x,\overline{3})\geq3.\;$Hence $x\preceq2$
and $C_{l-1}\preceq D_{l}$ by using that $C_{l-1}\preceq C_{l}.$ If $i_{1}=1$
and $h(C_{l})=2$, only the following configurations can appear%
\[
\left\{
\begin{tabular}
[c]{l}%
$\mathrm{(i):}$ $C_{l-1}C_{l}=%
\begin{tabular}
[c]{|l|l|}\hline
$x$ & $\mathtt{2}$\\\hline
$y$ & $\mathtt{3}$\\\hline
\end{tabular}
\text{ and }C_{l-1}D_{l}=%
\begin{tabular}
[c]{|l|l|}\hline
$x$ & $\mathtt{1}$\\\hline
$y$ & $\mathtt{3}$\\\hline
\end{tabular}
\vspace{0.15cm}$\\
$\mathrm{(ii):}$ $C_{l-1}C_{l}=%
\begin{tabular}
[c]{|l|l|}\hline
$x$ & $\mathtt{2}$\\\hline
$y$ & $\mathtt{0}$\\\hline
\end{tabular}
\text{ and }C_{l-1}D_{l}=%
\begin{tabular}
[c]{|l|l|}\hline
$x$ & $\mathtt{1}$\\\hline
$y$ & $\mathtt{3}$\\\hline
\end{tabular}
\vspace{0.15cm}$\\
$\mathrm{(iii):}$ $C_{l-1}C_{l}=%
\begin{tabular}
[c]{|l|l|}\hline
$x$ & $\mathtt{0}$\\\hline
$y$ & $\mathtt{0}$\\\hline
\end{tabular}
\text{ and }C_{l-1}D_{l}=%
\begin{tabular}
[c]{|l|l|}\hline
$x$ & $\mathtt{3}$\\\hline
$y$ & $\mathtt{0}$\\\hline
\end{tabular}
\vspace{0.15cm}$\\
$\mathrm{(iv):}$ $C_{l-1}C_{l}=%
\begin{tabular}
[c]{|l|l|}\hline
$x$ & $\mathtt{0}$\\\hline
$y$ & $\mathtt{\bar{2}}$\\\hline
\end{tabular}
\text{ and }C_{l-1}D_{l}=%
\begin{tabular}
[c]{|l|l|}\hline
$x$ & $\mathtt{3}$\\\hline
$y$ & $\mathtt{\bar{2}}$\\\hline
\end{tabular}
\vspace{0.15cm}$\\
$\mathrm{(v):}$ $C_{l-1}C_{l}=%
\begin{tabular}
[c]{|l|l|}\hline
$x$ & $\mathtt{\bar{3}}$\\\hline
$y$ & $\mathtt{\bar{2}}$\\\hline
\end{tabular}
\text{ and }C_{l-1}D_{l}=%
\begin{tabular}
[c]{|l|l|}\hline
$x$ & $\mathtt{3}$\\\hline
$y$ & $\mathtt{\bar{2}}$\\\hline
\end{tabular}
$%
\end{tabular}
\right.  .
\]
In case $\mathrm{(i)}$, and $\mathrm{(ii)}$ we have $x=1$ and $y\preceq3$ by
(\ref{cond_tab}).\ In cases $\mathrm{(iii),}$ $\mathrm{(iv)}$ and
$\mathrm{(v)}$ we have similarly $x\preceq3.\;$By using that $C_{l-1}\preceq
C_{l},$ we obtain $C_{l-1}\preceq D_{l}.$ If $i_{1}=1$ and $h(C_{l})=1$ only
the following configurations may appear:%
\[
\left\{
\begin{tabular}
[c]{l}%
$\mathrm{(vi):}C_{l-1}C_{l}=%
\begin{tabular}
[c]{|l|l}\hline
$x$ & \multicolumn{1}{|l|}{$\mathtt{0}$}\\\hline
$y$ & \\\cline{1-1}%
\end{tabular}
$ and $C_{l-1}D_{l}=%
\begin{tabular}
[c]{|l|l}\hline
$x$ & \multicolumn{1}{|l|}{$\mathtt{3}$}\\\hline
$y$ & \\\cline{1-1}%
\end{tabular}
\vspace{0.15cm}$\\
$\mathrm{(vii):}C_{l-1}C_{l}=$%
\begin{tabular}
[c]{|l|l|}\hline
$x$ & $\mathtt{0}$\\\hline
\end{tabular}
and $C_{l-1}D_{l}=$%
\begin{tabular}
[c]{|l|l|}\hline
$x$ & $\mathtt{0}$\\\hline
\end{tabular}
\end{tabular}
\right.
\]
with $x\preceq3$ by (\ref{cond_tab}). So $C_{l-1}\preceq D_{l}$ and the
proposition holds.
\end{proof}

Once the tableau $T_{1}$ defined, we do the same with $T_{1}$ getting a new
tableau $T_{2}$ of type $G_{2}$ and a new integer $i_{2}$. And so on until the
tableau $T_{s}$ obtained is equal to $T_{\lambda}$. Notice that we can not
write \textrm{w(}$T_{1})=\widetilde{e}_{i_{1}}^{r_{1}}\mathrm{w(}T)$ in
general, that is, our algorithm does not provide a path in the crystal graph
$B(\lambda)$ joining the vertex $\mathrm{w(}T)$ to the vertex of highest
weight $\mathrm{w(}T_{\lambda})$.

\begin{example}
Consider the tableau of type $G_{2}$%
\[
T=%
\begin{tabular}
[c]{|l|l|l}\hline
$\mathtt{3}$ & $\mathtt{\bar{3}}$ & \multicolumn{1}{|l|}{$\mathtt{\bar{1}}$%
}\\\hline
$\mathtt{\bar{2}}$ & $\mathtt{\bar{1}}$ & \\\cline{1-2}%
\end{tabular}
.
\]
We obtain successively:%
\[
T_{1}=%
\begin{tabular}
[c]{|l|l|l}\hline
$\mathtt{3}$ & $\mathtt{3}$ & \multicolumn{1}{|l|}{$\mathtt{\bar{2}}$}\\\hline
$\mathtt{\bar{2}}$ & $\mathtt{\bar{2}}$ & \\\cline{1-2}%
\end{tabular}
\text{, }T_{2}=%
\begin{tabular}
[c]{|l|l|l}\hline
$\mathtt{2}$ & $\mathtt{2}$ & \multicolumn{1}{|l|}{$\mathtt{\bar{3}}$}\\\hline
$\mathtt{\bar{3}}$ & $\mathtt{\bar{3}}$ & \\\cline{1-2}%
\end{tabular}
\text{, }T_{3}=%
\begin{tabular}
[c]{|l|l|l}\hline
$\mathtt{1}$ & $\mathtt{1}$ & \multicolumn{1}{|l|}{$\mathtt{3}$}\\\hline
$\mathtt{3}$ & $\mathtt{3}$ & \\\cline{1-2}%
\end{tabular}
\text{, }T_{4}=%
\begin{tabular}
[c]{|l|l|l}\hline
$\mathtt{1}$ & $\mathtt{1}$ & \multicolumn{1}{|l|}{$\mathtt{2}$}\\\hline
$\mathtt{2}$ & $\mathtt{2}$ & \\\cline{1-2}%
\end{tabular}
\]
and $A(T)=f_{1}^{(4)}f_{2}^{(5)}f_{1}^{(8)}f_{2}^{(3)}v_{T_{\lambda}}.$
\end{example}

\begin{proposition}
\label{prop_A(T)}The expansion of $A(T)$ on the basis $\{v_{\tau};\tau
\in\mathbf{T}(\lambda)\}$ of $W(\lambda)$ is of the form
\[
A(T)=\underset{\tau}{\sum}\alpha_{\tau,T}(q)v_{\tau}%
\]

where the coefficients $\alpha_{\tau,T}(q)$ satisfy:

\textrm{(i)}: $\alpha_{\tau,T}(q)\neq0$ only if $\tau$ and $T$ have the same weight,

\textrm{(ii)}: $\alpha_{\tau,T}(q)\in\mathbb{Z}[q,q^{-1}]$ and $\alpha_{T,T}(q)=1,$

\textrm{(iii)}: $\alpha_{\tau,T}(q)\neq0$ only if $\tau\trianglelefteq T.$
\end{proposition}

\begin{proof}
The proof is the same than in Proposition 4.3.4 of \cite{lec3}.
\end{proof}

It follows from \textrm{(iii)} that the vectors $A(T)$ are linearly
independent in $V(\lambda)$. This implies that $\{A(T);T\in\mathbf{T}%
_{G}(\lambda)\}$ is a $\mathbb{Q[}q]$-basis of $V(\lambda)$. Indeed by Theorem
\ref{KM}, $\dim V(\lambda)=\mathrm{card}(\mathbf{T}_{G}(\lambda))$. As a
consequence of (\ref{A(T)_monom}), we obtain $\overline{A(T)}=A(T)$.

\subsubsection{From $A(T)$ to $G(T)$}

To obtain $G(T)$ from $A(T)$ we proceed as in \cite{L-T}, \cite{L1} and
\cite{lec3}. The reader is referred to \cite{L-T} for the proofs.\ Set%
\[
G(T)=\underset{\tau\in\mathbf{T}_{G}(\lambda)}{\sum}d_{\tau,T}(q)\,v_{\tau}%
\]
Our aim is to describe a simple algorithm for computing the rectangular matrix
of coefficients
\[
D=[d_{\tau,T}(q)],\text{ \ \ \ \ }\tau\in\mathbf{T}(n,\lambda)\text{,
\ \ \ \ }T\in\mathbf{T}_{G}(\lambda).
\]

\begin{lemma}
\label{G_on_tabloids}The coefficients $d_{\tau,T}(q)$ belong to $\mathbb{Q}%
[q]$. Moreover $d_{\tau,T}(0)=0$ if $\tau\neq T$ and $d_{T,T}=1$.
\end{lemma}

Now write
\begin{equation}
G(T)=\underset{S\in T_{G}(\lambda)}{\sum}\beta_{S,T}(q)\,A(S)
\label{G(T)_on_A(T)}%
\end{equation}
the expansion of the basis $\{G(T)\}$ on the basis $\{A(T)\}$.

\begin{lemma}
The coefficients $\beta_{S,T}(q)$ of (\ref{G(T)_on_A(T)}) satisfy:

$\mathrm{(i)}$: $\beta_{S,T}(q)=\beta_{S,T}(q^{-1}),$

$\mathrm{(ii)}$: $\beta_{S,T}(q)=0$ unless $S\trianglelefteq T,$

$\mathrm{(iii)}$: $\beta_{T,T}(q)=1$.
\end{lemma}

Let $T_{\lambda}=T^{(1)}\vartriangleleft T^{(2)}\vartriangleleft\cdot
\cdot\cdot\vartriangleleft T^{(t)}$ be the sequence of tableaux of
$\mathbf{T}_{G}(\lambda)$ ordered in increasing order. We have $G(T_{\lambda
})=A(T_{\lambda})$, i.e. $G(T^{(1)})=A(T^{(1)})$. By the previous lemma, the
transition matrix $M$ from $\{A(T)\}$ to $\{G(T)\}$ is upper unitriangular
once the two bases are ordered with $\trianglelefteq$. Since $\{G(T)\}$ is a
$\mathbb{Q}[q,q^{-1}]$ basis of $V_{\mathbb{Q}}(\lambda)$ and $A(T)\in
V_{\mathbb{Q}}(\lambda),$ the entries of $M$ are in $\mathbb{Q}[q,q^{-1}].$
Suppose by induction that we have computed the expansion on the basis
$\{v_{\tau};\tau\in\mathbf{T}(\lambda)\}$ of the vectors%
\[
G(T^{(1)}),...,G(T^{(i)})
\]
and that this expansion verifies $d_{\tau,T^{(p)}}(q)=0$ if $\tau
\vartriangleright T^{(p)}$ for $p=1,...,i$. The inverse matrix $M^{-1}$ is
also upper unitriangular with entries in $\mathbb{Q}[q,q^{-1}]$. So we can
write:%
\begin{equation}
G(T^{(i+1)})=A(T^{(i+1)})-\gamma_{i}(q)G(T^{(i)})-\cdot\cdot\cdot-\gamma
_{1}(q)G(T^{(1)})\text{.} \label{G_on_A}%
\end{equation}
It follows from condition (\ref{cond_invo}) and Proposition \ref{prop_A(T)}
that $\gamma_{m}(q)=\gamma_{m}(q^{-1})$ for $m=1,...,i$. By Lemma
\ref{G_on_tabloids}, the coordinate $d_{T^{(i)},T^{(i+1)}}(q)$ of
$G(T^{(i+1)})$ on the vector $v_{T^{(i)}}$ belongs to $\mathbb{Q}[q]$,
$d_{T^{(i)},T^{(i+1)}}(0)=0$ and the coordinate $d_{T^{(i)},T^{(i)}}(q)$ of
$G(T^{(i)})$ on the vector $v_{T^{(i)}}$ is equal to $1$. Moreover
$v_{T^{(i)}}$ can only occur in $A(T^{(i+1)})-\gamma_{i}(q)G(T^{(i)}).$ If
\[
\alpha_{T^{(i)},T^{(i+1)}}(q)=\underset{j=-r}{\overset{s}{\sum}}a_{j}q^{j}%
\in\mathbb{Z}[q,q^{-1}]
\]
then we will have
\[
\gamma_{i}(q)=\overset{0}{\underset{j=-r}{\sum}}a_{j}q^{j}+\underset
{j=1}{\overset{r}{\sum}}a_{-j}q^{j}\in\mathbb{Z}[q,q^{-1}].
\]
Next if the coefficient of $v_{T^{(i-1)}}$ in $A(T^{(i+1)})-\gamma
_{i}(q)G(T^{(i)})$ is equal to
\[
\underset{j=-l}{\overset{k}{\sum}}b_{j}q^{j}%
\]
using similar arguments we obtain%
\[
\gamma_{i-1}(q)=\overset{0}{\underset{j=-l}{\sum}}b_{j}q^{j}+\overset
{l}{\underset{j=1}{\sum}}b_{-j}q^{j},
\]
and so on. So we have computed the expansion of $G(T^{(i+1)})$ on the basis
$\{v_{\tau}\}$ and this expansion verifies $d_{\tau,T^{(i+1)}}(q)=0$ if
$\tau\vartriangleright T^{(i+1)}$. Finally notice that $\gamma_{s}%
(q)\in\mathbb{Z}[q,q^{-1}]$ by Proposition \ref{prop_A(T)}.

\begin{theorem}
Let $T\in\mathbf{T}_{G}(\lambda).\;$Then $G(T)=\sum d_{\tau,T}(q)v_{\tau}$
where the coefficients $d_{\tau,T}(q)$ verify:

\textrm{(i)}: $d_{\tau,T}(q)\in\mathbb{Z}[q],$

\textrm{(ii)}: $d_{T,T}(q)=1$ and $d_{\tau,T}(0)=0$ for $\tau\neq T,$

\textrm{(iii)}: $d_{\tau,T}(q)\neq0$ only if $\tau$ and $T$ have the same
weight, and $\tau\trianglelefteq T$.
\end{theorem}

\end{document}